\documentclass[10pt]{amsart}

\usepackage{amssymb}
\usepackage{amsmath}
\usepackage{latexsym}

\theoremstyle{plain}
\newtheorem{theorem}{Theorem}[section]
\theoremstyle{plain}

\newtheorem{lemma}[theorem]{Lemma}
\newtheorem{proposition}[theorem]{Proposition}
\newtheorem{corollary}[theorem]{Corollary}
\theoremstyle{definition}

\theoremstyle{remark}
\newtheorem{remark}{Remark}[section]

\newcommand{\dt}{\Delta t}
\newcommand{\dx}{\Delta x}
\newcommand{\A}{\ensuremath{\mathcal{A}}}
\newcommand{\B}{\ensuremath{\mathcal{B}}}
\newcommand{\eps}{\ensuremath{\varepsilon}}
\newcommand{\R}{\ensuremath{\mathbb{R}}}
\newcommand{\N}{\ensuremath{\mathbb{N}}}
\newcommand{\Z}{\ensuremath{\mathbb{Z}}}
\newcommand{\CS}{\mathcal{S}}
\newcommand{\BS}{\ensuremath{\mathbb{S}}}
\newcommand{\del}{\partial}
\newcommand{\ol}{\overline}
\newcommand{\ra}{\rightarrow}
\newcommand{\alp}{\alpha}
\newcommand{\LL}{\mathcal{L}}
\newcommand{\II}{\mathcal{I}}
\newcommand{\M}{\mathcal{M}}
\newcommand{\G}{\mathcal G}
\newcommand{\tr}{\mathrm{tr}}
\newcommand{\bx}{\bar{x}}
\newcommand{\by}{\bar{y}}
\newcommand{\bt}{\bar{t}}
\newcommand{\bu}{\bar{u}}
\newcommand{\bv}{\bar{v}}
\newcommand{\ba}{\bar{a}}
\newcommand{\bb}{\bar{b}}
\newcommand{\bc}{\bar{c}}
\newcommand{\bff}{\bar{f}}
\newcommand{\bs}{\bar{\sigma}}
\newcommand{\CC}{\mathcal{C}}

\newcommand{\argmin}{\mathrm{argmin}}

\def\sym{\BS^N}
\def\e{\varepsilon}
\def\tb{\bar t}
\def\xb{\bar x}

\numberwithin{equation}{section}
\allowdisplaybreaks

\title[Error bounds]{Error bounds for monotone approximation schemes
  for parabolic Hamilton-Jacobi-Bellman equations}

\date{\today}

\author[Barles]{Guy Barles}
\address[Guy Barles]{\newline
Laboratoire de Math\'ematiques et Physique Th\'eorique\newline
University of Tours\newline
Parc de Grandmont \newline
37200 TOURS, France}
\email[]{barles\@@lmpt.univ-tours.fr}
\urladdr{http://www.phys.univ-tours.fr/\~{}barles}

\author[Jakobsen]{Espen R.~Jakobsen}
\address[Espen R.~Jakobsen]{\newline
    Department of Mathematical Sciences\newline
    Norwegian University of Science and Technology\newline
    N--7491 Trondheim, Norway}
\email[]{erj\@@math.ntnu.no}
\urladdr{http://www.math.ntnu.no/\~{}erj}


\thanks{Jakobsen was supported by the Research
Council of Norway, grant no. 151608/432.}

\begin{document}

\begin{abstract} We obtain non-symmetric upper and lower bounds on the
  rate of convergence of general monotone approximation/numerical
  schemes for parabolic Hamilton Jacobi Bellman Equations by
  introducing a new notion of consistency. We apply our general
  results to various schemes including finite difference schemes,
  splitting methods and the classical approximation by 
  piecewise constant controls. 
\end{abstract}

\maketitle

\section{Introduction}
\label{Sec:intro}

In this article, we are interested in the rate of convergence of
general monotone approximation/numerical schemes for time-dependent
Hamilton Jacobi Bellman (HJB) Equations.

In order to be more specific, the HJB Equations we consider are written in the
following form
\begin{align}
\label{E}
u_t+F(t,x,u,Du,D^2u)&=0&&\text{in}\quad Q_T:=(0,T]\times\R^N,\\
u(0,x)&=u_0(x)&&\text{in}\quad\R^N,
\label{IV}
\end{align}
where
$$ F(t,x,r,p,X)=\sup_{\alpha\in \A}\left\{\LL^{\alp} (t,x,r,p,X)\right\}\; ,$$
with
$$ \LL^{\alp} (t,x,r,p,X) := - \tr[a^{\alpha}(t,x) X] -
b^{\alpha}(t,x) p - c^{\alpha}(t,x) r - f^{\alpha}(t,x).$$ 
The coefficients $a^{\alpha}$, $b^{\alpha}$,
$c^{\alpha}$, $f^{\alpha}$ and the initial data $u_0$ take values respectively in
$\sym$, the space of $N \times N$ symmetric matrices, $\R^N$, $\R$,
$\R$, and $\R$. Under suitable assumptions (see  {\bf (A1)} in
Section~\ref{Sec:Sw}),  
the initial value problem \eqref{E}-\eqref{IV} has a unique, bounded,
H\"older continuous, viscosity solution $u$ which is the value
function of a finite horizon, optimal stochastic control problem.

We consider approximation/numerical schemes for \eqref{E}-\eqref{IV}
written in the following abstract way
\begin{align}
\label{S}
S(h,t,x,u_h(t,x),[u_h]_{t,x})&=0 &&\text{in}\quad
\G_h^+:=\G_h\setminus\{t=0\},\\ 
u_h(0,x)&=u_{h,0}(x) && \text{in}\quad \G_h^0:=\G_h\cap\{t=0\},\nonumber
\end{align}
where $S$ is, loosely speaking, a consistent, monotone and uniformly
continuous approximation of the equation \eqref{E} defined on a
grid/mesh $\G_h\subset\ol Q_T$.
The approximation parameter $h$ can be multi-dimensional,
e.g. $h$ could be $(\dt,\dx)$, $\dt$, $\dx$ denoting time and space
discretization parameters, $\dx$ can be itself multi-dimensional. 
The approximate solution is $u_h:\G_h\ra \R$, $[u_h]_{t,x}$ is a
function defined from $u_h$ representing, typically, the value of
$u_h$ at other points than $(t,x)$. We assume that the total scheme
including the initial value is well-defined on some appropriate subset
of the space of bounded continuous functions on $\G_h$.

The abstract notation was introduced by Barles and Souganidis
\cite{BS:Conv} to display clearly the monotonicity of the scheme. One
of the main assumptions is that $S$ is non-decreasing in $u_h$ and
non-increasing in $[u_h]_{t,x}$ with the classical ordering of
functions. The typical approximation schemes we have in mind are
various finite differences numerical scheme (see e.g. Kushner
and Dupuis \cite{KD:Book} and Bonnans and Zidani \cite{BZ}) and
control schemes based on the dynamic programming principle (see e.g.
Camilli and Falcone \cite{CF:Appr}). However, for reasons 
explained below, we will not discuss control schemes in this paper.

The aim of this paper is to obtain estimates on the rate of the
convergence of $u_h$ to $u$. To obtain such results, one faces the
double difficulty of having to deal with both  \emph{fully nonlinear}
equations and {\em non-smooth} solutions. Since these
equations may be also degenerate, the (viscosity) solutions are
expected to be no more than H\"older continuous in general.

Despite of these difficulties, in the 80's,  Crandall \& Lions
\cite{CL:Appr} provided the first optimal rates of convergence for first-order
equations. We refer to Souganidis \cite{So:Appr} for more general
results in this direction. For technical reasons, the problem
turns out to be more difficult for second-order equations, and the
question remained open for a long time. 

The breakthrough came in 1997 and 2000 with Krylov's papers
\cite{Kr:HJB1,Kr:HJB2}, and by now there exists several papers based
on and extending his ideas,
e.g. \cite{BJ:Rate,BJ:Rate2,DK:ConstCoeff,J:Par,Kr:Const,Kr:LipCoeff}.
The main idea of Krylov is a method named by himself ``shaking the
coefficients''. Combined with a standard mollification
argument, it allows one to get smooth {\em subsolutions} of the
equation which approximate the solution. Then classical arguments
involving consistency and monotonicity of the scheme yield a {\em
  one-sided bound} on the error. This method uses in a crucial way the
convexity of the equation in $u$, $Du$, and $D^2 u$.

It is much more difficult to obtain the other bound and essentially there
are two main approaches. The first one consists of interchanging the
role of the scheme and the equation. By applying the above explained
ideas, one gets a sequence of appropriate smooth subsolutions of the scheme and
concludes by consistency and the comparison principle for the
equation. This idea was used in different articles, see
\cite{BJ:Rate,DK:ConstCoeff,J:Par,Kr:HJB1,Kr:LipCoeff}. Here, the key
difficulty is to obtain a ``continuous dependence'' result
for the scheme. Even though it is now standard to prove that the solutions
of the HJB Equation with ``shaken coefficients'' remain close to the
solution of the original equation, such type of results are not known
for numerical schemes in general. We mention here the nice paper of 
Krylov \cite{Kr:LipCoeff} where such kind of results are obtained by a
tricky Bernstein type of argument.
However, these results along with the corresponding error bounds, only
hold for equations and schemes with special structures.

The second approach consists of considering some approximation of the
equation or the associated control problem and to obtain the other bound
either by probabilistic arguments (as Krylov first did using piecewise
constant controls, \cite{Kr:Const,Kr:HJB2}) or by 
building a sequence of appropriate ``smooth supersolution'' of the
equation (see \cite{BJ:Rate2} where, as in the present paper,
approximations by switching are considered).

The first approach leads to better error bounds than the
second one but it seems to work only for very specific schemes and
with restrictions on the equations. The second approach yields
error bounds in ``the general case'' but at the expense of lower
rates. 

In this paper we use the second approach by extending the methods
introduced in \cite{BJ:Rate2}. Compared with the various results of
Krylov, we obtain better rates in most cases, our results apply to
more general schemes, and we use a simpler, purely analytical approach. In 
fact our method is robust in the sense that it applies to ``general''
schemes without any particular form and under rather natural
assumptions. However, we mention again that in certain situations
the first approach can be used to get better rates, see in particular
\cite{Kr:LipCoeff}.

The results in \cite{BJ:Rate2} apply to stationary HJB equations set
in whole space $\R^N$. In this paper we extend these results to
initial value problems for time-dependent HJB equations. The latter
case is much more interesting in view of applications, and from a
mathematical point of view, slightly more difficult. However, in our opinion
the most important difference between the two papers lays in the
formulation of the consistency requirements 
and the main (abstract) results. Here we introduce a new (and more
general) formulation that emphasizes more the non-symmetrical feature of the
upper and lower bounds and their proofs. It is a kind of a recipe on
how to obtain error bounds in different situations, one which we feel
is easier to apply to new problems and gives better insight into how
the error bounds are produced. We also present several technical
improvements and simplifications in the proofs and, finally, several
new applications, some for which error bounds have not appeared
before: Finite difference methods 
(FDMs) using the $\theta$-method for time discretization, semidiscrete
splitting methods, and approximation by piecewise constant controls.

The results for finite difference approximations can be compared with
the ones obtained by Krylov in \cite{Kr:HJB2,Kr:Const}. As in
\cite{BJ:Rate2}, we get the rate $1/5$ for monotone FDMs while the
corresponding result in \cite{Kr:Const} is $1/21$. Of course, in
special situations the rate can be improved to $1/2$ which is the
optimal rate under our assumptions. We refer to \cite{Kr:LipCoeff} for
the most general results in that direction, and to \cite{DK:Lin} for
the optimality of the rate $1/2$. The results for semidiscrete
splitting methods are new, while the ones for the control
approximation we get $1/10$ which is worse than $1/6$ obtained by
Krylov in \cite{Kr:Const}. It would be interesting to understand why
Krylov is doing better than us here but not in the other cases.

We conclude this introduction by explaining the notations we will use
throughout this paper.  By $|\cdot|$ we mean the standard Euclidean
norm in any $\R^p$ type space (including the space of $N \times P$
matrices). In particular, if $X\in \sym$, then $|X|^2=\tr(XX^T)$ where
$X^T$ denotes the transpose of $X$.

If $w$ is a bounded function from some set $Q' \subset \ol Q_\infty$
into either $\R$, $\R^M$, or 
the space of $N \times P$ matrices, we set
$$ |w|_{0} = \sup_{(t,y)\in Q'}\,|w(t,y)|.$$ Furthermore, for
$\delta\in(0,1]$, we set
$$[w]_{\delta} = \sup_{(t,x)\neq
  (s,y)}\,\frac{|w(t,x)-w(s,y)|}{(|x-y|+|t-s|^{1/2})^\delta}\quad
\hbox{and}\quad |w|_{\delta}=  |w|_{0} + [w]_{\delta}.$$
Let $C_b (Q')$ and $\CC^{0,\delta}(Q')$, $\delta\in(0,1],$ denote 
  respectively the space of 
bounded continuous functions on $Q'$ and the subset of $C_b (Q')$ in
which the norm  $|\cdot|_\delta$ is finite.
Note in particular the choices $Q'=Q_T$ and $Q'=\R^N$. In the
following we always suppress the domain $Q'$ when writing norms.

We denote by $\leq$ the component by component ordering in $\R^M$ and
the ordering in the sense of positive semi-definite matrices in $\sym$.
For the rest of this paper we let $\rho$ denotes the same, fixed,
positive smooth function with support in $\{0<t<1\}\times\{|x|<1\}$ and mass
$1$. From this function $\rho$, we define the sequence of mollifiers 
$\{\rho_\eps\}_{\eps>0}$ as follows,
$$\rho_\eps(t,x)=\frac{1}{\eps^{N+2}}\rho\left(\frac{t}{\eps^2},
\frac{x}{\eps}\right) \quad \text{in}\quad \ol Q_\infty.$$

The rest of this paper is organized as follows: In the next section we present
results on the so-called switching approximation for the problem
\eqref{E}-\eqref{IV}. As in \cite{BJ:Rate2}, these results are crucial
to obtain the general results on the rate of convergence of
approximation/numerical schemes and are of an independent
interest. Section~\ref{Sec:HJB} is devoted to state and prove the main
result on the rate of convergence. Finally we present various
applications to classical finite difference schemes, splitting method
and on the classical approximation by piecewise constant controls.

\section{Convergence Rate for a Switching System}
\label{Sec:Sw}

In this section, we obtain the rate of convergence for a certain switching
system approximations to the HJB equation \eqref{E}. Such
approximations have be studied in \cite{EF:Sw,CDE:Sw}, and a viscosity
solutions theory of switching systems can be found in
\cite{Ya:Sw2,IK:Sw,IK:Sys}. We consider the following type of
switching systems, 
\begin{align}
\label{Sw1}
F_i(t,x,v,\del_tv_i,Dv_i,D^2v_i)&=0 && \text{in}\quad Q_T, \quad
i\in\II:=\{1,\dots,M\},\\
v(0,x)&=v_0(x)&&\text{in}\quad\R^N,\nonumber
\end{align}
where the solution $v=(v_{1}, \cdots ,v_{M})$ is in $\R^M$, and for 
$i\in\II$, 
$(t,x)\in Q_T$, $r =(r_{1}, \cdots ,r_{M}) \in \R^M$, $p_t\in\R$, $p_x\in
\R^N$, and $X \in \mathcal{S}^N$, $F_{i}$ is given by
$$
F_i(t,x,r,p_t,p_x,X)
=\max\Big\{p_t+\sup_{\alp\in\A_{i}}\LL^{\alp}(t,x,r_i,p_x,X);
r_i-\M_ir\Big\},$$
where the $\A_{i}$'s are subsets of $\A$, $\LL^{\alp}$ is defined below \eqref{E}, and
for $k>0$, 
$$
\M_ir =\min_{j\neq i}\{r_j+k\}.
$$ 
Finally for the initial data, we are interested here in the case
when $v_0=(u_0,\dots,u_0)$.

Under suitable assumptions on the data (See {\bf (A1)} below), we have
existence and uniqueness of a solution $v$ of this system.  Moreover,
it is not so difficult to see that, as $k\to 0$, every component of
$v$ converge locally uniformly to the solution of the following HJB
equation
\begin{align}
\label{HJBi}
u_t+\sup_{\alp\in \tilde \A} \LL^{\alp}(x,u,Du,D^2u)&=0 && \text{in}
\quad Q_T,\\
u(0,x)&=u_0(x) && \text{in}\quad \R^N,\nonumber
\end{align}
where $ \tilde \A = \cup_{i}\, \A_{i}$.  

The objective of this section is to obtain an error bound for this
convergence. For the sake of simplicity, we restrict ourselves to the
situation where the solutions are in $\CC^{0,1}(Q_T)$,
i.e. when they are bounded, Lipschitz continuous in $x$, and H{\"o}lder
$1/2$ in $t$. Such type of regularity is natural in this context.
 However, it is not difficult to adapt our approach to more
general situations, and we give results in this direction in
Section~\ref{Sec:Rem}. 

We will use the following assumption
\medskip

\noindent {\bf (A1)} For any $\alp \in \A$,
$a^{\alp}=\frac12\sigma^{\alp}{\sigma^{\alp}}^T$ for some $N\times P$ 
matrix $\sigma^{\alp}$. Moreover, there is a constant $K$ independent
of $\alp$ such that
$$|u_0|_1+|\sigma^{\alp}|_1+|b^{\alp}|_1 + |c^{\alp}|_1+|f^{\alp}|_1 \leq
K.$$

Assumption (A1) ensures the well-posedness of all the equations and
systems of equations we consider in this paper; we refer the reader to
the Appendix for a (partial) proof of this claim. In the present
situation, we have the following well-posedness and regularity result.
\begin{proposition}
\label{WPi}
Assume (A1). Then there exist unique solutions $v$ and
$u$ of \eqref{Sw1} and \eqref{HJBi} respectively, satisfying
$$|v|_1+|u|_1\leq C,$$
where the constant $C$ only depends on $T$ and $K$ appearing in (A1). 

Furthermore, if $w_1$ and $w_2$ are sub- and supersolutions of
\eqref{Sw1} or \eqref{HJBi} satisfying $w_1(0,\cdot)\leq w_2(0,\cdot)$,
then $w_1\leq w_2$.
\end{proposition}
\begin{remark}
\label{extension}
The functions $\sigma^{\alp}, b^{\alp}, c^{\alp}, f^{\alp}$ are a
priori only defined for times $t \in [0,T]$. But they can easily be extended
to times $[-r,T+r]$ for any $r\in\R^+$ in such a way that (A1) still holds.
In view of Proposition \ref{WPi} we can then solve our initial value problems
\eqref{Sw1} and \eqref{HJBi} either up to time $T+r$ and even, by
using a translation in time, on time intervals of the form
$[-r,T+r]$. We will use this fact several times below. 
\end{remark}

In order to obtain the rate of convergence for the switching
approximation, we use a regularization procedure introduced by Krylov
\cite{Kr:HJB2,BJ:Rate}.  This procedure requires the following
auxiliary system 
\begin{align}
\label{Sw2}
F^\eps_i(t,x,v^\eps,\del_t v_i^\eps,Dv^\eps_i,D^2v^\eps_i)&=0 
 &&\text{in}\quad Q_{T+\eps^2}, \quad i\in\II,\\
v^\eps(0,x) &= v_0(x) && \text{in} \quad \R^N,\nonumber
\end{align}
where $v^\eps = (v^\eps_1, \cdots, v^\eps_M)$,
\begin{align*}
&F^\eps_i(t,x,r,p_t,p_x,M)=\\
&\max\Big\{p_t+\sup_{
{\displaystyle{\mathop{\scriptstyle{\alp\in\A_{i}}}_{0 \leq s  \leq\eps^2, |e|\leq\eps}}}
 }\LL^{\alp}(t+s,x+e,r_i,p_x,X);\ 
r_i-\M_ir\Big\},  
\end{align*}
and $\LL$ and $\M$ are defined below \eqref{E} and \eqref{Sw1}
respectively. Note that we use here the extension mentioned in Remark
\ref{extension}.

By Theorems~\ref{WP} and \ref{CD} in the Appendix, we
have the following result:
\begin{proposition}
\label{WPi2}
Assume (A1). Then there exist a unique solution $v^\eps:\ol Q_{T+\eps^2}\ra \R$
of \eqref{Sw2} satisfying 
$$|v^\eps|_1+\frac1\eps|v^\eps-v|_0 \leq C,$$
where $v$ solves \eqref{Sw1} and the constant $C$ only depends on $T$ and
$K$ from (A1). 

Furthermore, if $w_1$ and $w_2$ are sub- and
supersolutions of \eqref{Sw2} satisfying $w_1(0,\cdot)\leq
w_2(0,\cdot)$, then $w_1\leq w_2$. 
\end{proposition}

We are now in a position to state and prove the main result of this section.
\begin{theorem}
\label{sw-rate}
Assume (A1) and $v_0=(u_0,\dots,u_0)$. If $u$ and $v$ are the solutions
of \eqref{HJBi} and \eqref{Sw1} respectively, then for $k$ small enough, 
$$0\leq v_i - u\leq Ck^{1/3}\quad \text{in}\quad Q_T,\quad i\in\II,$$
where $C$ only depends on $T$ and $K$ from (A1).
\end{theorem}

\begin{proof}
Since $w=(u,\dots,u)$ is a subsolution of \eqref{Sw1},
comparison for \eqref{Sw1} (Proposition~\ref{WPi}) yields $u\leq
v_i$ for $i\in\II$. 

To get the other bound, we use an argument suggested by P.-L. Lions
\cite{Li:PC} together with the regularization procedure of Krylov
\cite{Kr:HJB2}.
Consider first system \eqref{Sw2}. It follows that, for every $0 \leq
s  \leq\eps^2, |e|\leq\eps$, 
\begin{align*}
\del_t v^\eps_i+\sup_{\alp\in\A_{i}}\LL^\alp
(t+s,x+e,v^{\eps}_i(t,x),Dv^{\eps}_i,D^2v^{\eps}_i)\leq 
0 \quad \text{in}\quad Q_{T+\eps^2}, \quad i\in\II. 
\end{align*}
After a change of variables, we see that for every $0 \leq s
\leq\eps^2, |e|\leq\eps$, $v^{\eps}(t-s,x-e)$ is a subsolution of the
following system of uncoupled equations
\begin{align}
\label{LSys}
\del_tw_i+\sup_{\alp\in\A_{i}}\LL^\alp (t,x,w_i,Dw_i,D^2w_i)= 0 \quad
\text{in}\quad  Q^\eps_T, \quad i\in\II,
\end{align}
where $Q^\eps_T:=(\eps^2,T)\times\R^N$.
Define ${v_\eps}:=v^\eps*\rho_\eps$ where $\{\rho_\eps\}_{\eps}$ is the
sequence of mollifiers defined at the end of the introduction.  
A Riemann-sum approximation shows that $v_\eps(t,x)$
can be viewed as the limit of convex combinations of
$v^{\eps}(t-s,x-e)$'s for $0<s<\eps^2$ and $|e|<\eps$. Since the
$v^{\eps}(t-s,x-e)$'s 
are subsolutions of the \emph{convex} equation \eqref{LSys}, so are the
convex combinations.  By the stability result for viscosity subsolutions 
we can now conclude that $v_{\eps}$ is itself a
subsolution of \eqref{LSys}. 
We refer to the Appendix in \cite{BJ:Rate} for more details.

On the other hand, since $v^\eps$ is a continuous subsolution of 
\eqref{Sw2}, we have
$$ v^{\eps}_i  \leq \min_{j\neq i}v^{\eps}_j  +
k\quad\text{in}\quad Q_{T+\eps^2},\quad i\in\II.$$ 
It follows that ${\rm max}_{i}\,v^{\eps}_i (t,x) - {\rm
  min}_{i}\,v^{\eps}_i (t,x) \leq k$ in $Q_{T+\eps^2}$, and hence
$$ |v^{\eps}_i-v^{\eps}_j|_0\leq k, \quad i,j\in\II.$$
Then, by the definition and properties of $v_{\eps}$, we have
\begin{align*}
|\del_t{v_\eps}_{i}-\del_t{v_\eps}_{j}|_0\leq C\frac k {\eps^2},\quad
|D^n{v_\eps}_{i}-D^n{v_\eps}_{j}|_0\leq C\frac k {\eps^n}, \quad
n\in\N,\quad i,j\in\II,
\end{align*}
where $C$ depends only on $\rho$ and the uniform bounds on
${v_\eps}_{i}$ and $D{v_\eps}_{i}$, i.e. on $T$ and $K$ given in
(A1). Furthermore, from these bounds, we see that for $\eps<1$,
$$
\left|\del_t{v_\eps}_{j}+\sup_{\alp\in\A_{i}}\LL^\alp [{v_\eps}_j] -
\del_t{v_\eps}_{i}-\sup_{\alp\in\A_{i}}\LL^\alp[{v_\eps}_i]\right|
\leq  C\frac{k}{\eps^2} \quad 
\text{in}\quad Q^\eps_T,\quad i,j\in\II.$$
Here, as above, $C$ only depends on $\rho$, $T$ and $K$.
Since $v_{\eps}$ is a subsolution of \eqref{LSys}, this means that,
$$\del_t{v_\eps}_i+\sup_{\alp\in\ol\A}\LL^\alp(x,{v_\eps}_i,D{v_\eps}_i,
D^2{v_\eps}_i)  \leq C\frac{k}{\eps^2}\quad \text{in}\quad Q_T^\eps,\quad
i\in\II.$$ 
From assumption (A1) and the structure of the equation, we see that
${v_\eps}_i - t e^{Kt} C\frac{k}{\eps^2}$ is a 
subsolution of equation \eqref{HJBi} restricted to $Q_T^\eps$.

Comparison for \eqref{HJBi} restricted to $Q_T^\eps$
(Proposition~\ref{WPi}) yields 
$${v_\eps}_i-u\leq e^{Kt}\left(|{v_\eps}_i(\eps^2,\cdot)-u(\eps^2,\cdot)|_0+
Ct\frac{k}{\eps^2}\right)\quad \text{in}\quad Q_T^\eps, \quad i\in\II.$$
Regularity of $u$ and $v_i$ (Proposition \ref{WPi}) implies that
\begin{align*}
&|u(t,\cdot)-v_i(t,\cdot)|_0\leq ([u]_1+[v_i]_1) \eps \qquad \text{in}\quad
[0,\eps^2].
\end{align*}
Hence by Proposition \ref{WPi2}, regularity of $u$ and $v^\eps_i$,
and properties of mollifiers, we have  
$$v_i - u\leq v_i-{v_\eps}_i+{v_\eps}_i-u\leq C(\eps + \frac{k}{\eps^2}) \quad 
\text{in}\quad Q^\eps_T, \quad i\in\II.$$
Minimizing w.r.t. $\eps$ now yields
the result.
\end{proof}

\section{Convergence rate for the HJB equation}
\label{Sec:HJB}

In this section we derive our main result, an error bound for the
convergence of the solution of the scheme \eqref{S} to the solution of
the HJB Equation \eqref{E}-\eqref{IV}. 
As in \cite{BJ:Rate2}, this result is general and
derived using only PDE methods, and it extends and improves earlier
results by Krylov \cite{Kr:HJB1,Kr:HJB2}, Barles and Jakobsen
\cite{BJ:Rate,J:Par}. Compared to \cite{BJ:Rate2}, we consider here
the time-dependent case and introduce a new, improved, formulation of
the consistency requirement.

Throughout this section, we assume that (A1) holds and we
recall that, by Proposition \ref{WPi}, there exists a unique
$\CC^{0,1}$-solution $u$ of \eqref{E} satisfying $|u|_1\leq 
C$, where the constant $C$ only depends on $T$ and $K$ from (A1).
In Section \ref{Sec:Rem}, we will weaken assumption (A1) and give results for
$\CC^{0,\beta}$ solutions, $\beta\in(0,1)$.

In order to get a lower bound bound on the error, we have to require a
technical assumption: If $\{\alp_i\}_{i\in\II}\subset\A$ is a
sufficiently refined grid for $\A$, the solution associated to the
control set $\{\alp_i\}_{i\in\II}$ is close to $u$. In fact for this
to be true we need to assume that the coefficients
$\sigma^\alp,b^\alp,c^\alp,f^\alp$ can be approximated uniformly in
$(t,x)$ by $\sigma^{\alp_i},b^{\alp_i},c^{\alp_i},f^{\alp_i}$. The
precise assumption is:

\medskip
\noindent{\bf (A2)} For every $\delta>0$, there are $M\in\N$ and
$\{\alp_i\}_{i=1}^M\subset\A$, such that for any $\alp \in \A$,
$$ \inf_{1\leq i\leq
M}\left(|\sigma^\alp-\sigma^{\alp_i}|_{0}+|b^\alp-
b^{\alp_i}|_{0}+|c^\alp-c^{\alp_i}|_{0} 
+|f^\alp-f^{\alp_i}|_{0}\right)<\delta.$$\\[-0.5cm]

\smallskip
We point out that this assumptions is automatically satisfied if
either $\A$ is a finite set or if $\A$ is compact and $\sigma^\alp,
b^\alp, c^\alp, f^\alp $ are uniformly continuous functions of 
$t$, $x$, and $\alp$.

Next we introduce the following assumptions for the scheme \eqref{S}.

\smallskip
\noindent {\bf (S1)} {\bf (Monotonicity)}  There exists $\lambda,
\mu\geq 0 ,h_0>0$ such that if $|h| \leq h_0$, $u\leq v$ are functions
in $C_b(\G_h)$, and $\phi(t)=e^{\mu t}(a+bt)+c$ for
$a,b,c\geq0$, then  
$$S(h,t,x,r+\phi(t),[u+\phi]_{t,x}) \geq  S(h,t,x,r,[v]_{t,x}) + b/2 -
\lambda c
\quad \text{in} \quad \G_h^+.
$$

\medskip

\noindent {\bf (S2)} {\bf (Regularity)}
For every $h$ and $\phi\in C_b(\G_h)$, the function 
$(t,x) \mapsto$\linebreak $S(h,t,x,\phi(t,x),[\phi]_{t,x})$
is bounded and continuous in $\G_h^+$ and the function $r \mapsto
S(h,t,x,r,[\phi]_{t,x})$ is uniformly continuous for bounded $r$,
uniformly in $(t,x) \in \G_h^+$.
\begin{remark}
\label{notC_b}
In (S1) and (S2) we may replace $C_b(\G_h)$ by any relevant subset of
this space. The point is that \eqref{S} has to make sense for the
class of functions used. In Section \ref{Sec:FDM}, $C_b(\R^N)$ is itself
the relevant class of functions, while, in Section \ref{Sec:Appl}, it is
$C(\{0,1,\dots,n_T\};\CC^{0,1}(\R^N))$ (since $\G_h=\{0,1,\dots,n_T\}\times\R^N$). 
\end{remark}
Assumptions (S1) and (S2) imply a comparison result for the scheme
\eqref{S}, see Lemma \ref{compsche} below. 

Let us now state the key consistency conditions.
\smallskip

\noindent {\bf (S3)(i)} {\bf (Sub-consistency)} There exists
a function $E_1(\tilde K,h,\e)$ such
that for any sequence $\{\phi_\eps\}_{\eps>0}$ of smooth functions satisfying
$$  |\partial_t^{\beta_0}D^{\beta'} \phi_\e (x,t) | \leq \tilde K
\e^{1-2\beta_0-|\beta'|}  \quad \hbox{in  }\overline Q_T ,\quad\text{ for
any $\beta_0\in\N$, $\beta'=(\beta'_i)_i\in\N^{N}$},$$
where $|\beta'|=\sum_{i=1}^N \beta_i'$, the following inequality holds:
$$ S(h,t,x,\phi_\e (t,x),[\phi_\e ]_{t,x})Ê\leq \phi_{\e t}+
F(t,x,\phi,D\phi_\e,D^2\phi_\e) + E_1 (\tilde K,h,\e)\quad \hbox{in  }\G_h^+.$$

\smallskip
\noindent {\bf (S3)(ii)} {\bf (Super-consistency)} There exists
a function $E_2(\tilde K,h,\e)$ such
that for any sequence $\{\phi_\eps\}_\eps$ of smooth functions satisfying 
$$  |\partial_t^{\beta_0}D^{\beta'} \phi_\e (x,t) | \leq \tilde K
\e^{1-2\beta_0-|\beta'|}  \quad \hbox{in  }\overline Q_T ,\quad\text{ for
any $\beta_0\in\N$, $\beta'\in\N^{N}$},$$
the following inequality holds:
$$ S(h,t,x,\phi_\e (t,x),[\phi_\e ]_{t,x})Ê\geq \phi_{\e t} +
F(t,x,\phi,D\phi_\e,D^2\phi_\e) - E_2 (\tilde K,h,\e)\quad \hbox{in  }\G_h^+.$$

\medskip
Typically the $\phi_\eps$ we have in mind in (S3) are of the form $\chi_\eps *
\rho_\eps$ where $(\chi_\eps)_\eps$ is a sequence of uniformly bounded functions in $\CC^{0,1}$ and $\rho_\eps$ is the mollifier
defined at the end  of the introduction.

The main result in this paper is the following:
\begin{theorem}
\label{mainres}
Assume (A1), (S1), (S2) and that \eqref{S} has a unique
solution $u_h$ in $C_b(\G_h)$. Let $u$ denote the solution of
\eqref{E}-\eqref{IV}, and let $h$ be sufficiently small. 

(a) {\bf (Upper bound)} If (S3)(i) holds, then there exists a constant
$C$ depending only $\mu$, $K$ in (S1), (A1) such that
$$ u-u_h \leq e^{\mu t}|(u_0-u_{0,h})^+|_0 + C\min_{\e>0} \left(\e + E_1 (\tilde
K,h,\e)\right) \quad\text{in}\quad\G_h,
$$
where $\tilde K = |u|_1$.

(b) {\bf (Lower bound)} If (S3)(ii) and (A3) holds, then there exists
a constant $C$ depending only $\mu$, $K$ in (S1), (A1) such that
$$
u-u_h \geq -e^{\mu t}|(u_0-u_{0,h})^-|_0 - C\min_{\e>0} \left(\e^{1/3} + E_2
(\tilde K,h,\e)\right) \quad\text{in}\quad\G_h, 
$$
where $\tilde K = |u|_1$.
\end{theorem}

\medskip

The motivation for this new formulation of the upper and lower bounds
is threefold: (i) in some applications, $E_1\neq E_2$ and therefore it is natural to have such disymmetry
in the consistency requirement (see Section \ref{Sec:Appl}), (ii) from the proof it can be seen that the
upper bound (a) is proven independently of the lower bound (b), and
most importantly, (iii) the new formulation describes completely how
the bounds are obtained from the consistency requirements.
The good $h$-dependence and the bad $\e$ dependence of $E_1$ and $E_2$ are
combined in the minimization process to give the final bounds, see
Remark \ref{S3} below. 

Since the minimum is achieved for $\e \ll 1$, the upper bound is in
general much better than the lower bound (in particular in cases where
$E_1 = E_2$). 

Finally note that the existence of a $u_h$ in $C_b(\G_h)$
must be proved for each particular scheme $S$. We refer to
\cite{Kr:HJB1,Kr:HJB2,BJ:Rate,J:Par} for examples of such
arguments.

\begin{remark}
\label{S3}
In the case of a finite difference method with a time step $\dt$ and
maximal mesh size in space $\dx$, a standard formulation of the
consistency requirement would be 

\medskip

\noindent {\bf (S3')} There exist finite sets
$I\subset \N\times\N_0^N, \bar I \subset \N_0\times\N^{N}$ and
constants $K_c\geq0$, $k_\beta,\bar 
k_{\bar\beta}$ for $\beta=(\beta_0,\beta')\in I, \bar\beta=(\bar\beta_0,\bar\beta')\in \bar I$ such
that for every  $h=(\dt,\dx) > 0$, $(t,x)\in \G_h^+$, and smooth
functions $\phi$:
\begin{align*}
&\left|\phi_t +
F(t,x,\phi,D\phi,D^2\phi)-S(h,t,x,\phi(t,x),[\phi]_{t,x})\right|\\
&\leq K_c\sum_{\beta\in
  I}|\del_t^{\beta_0}D^{\beta'}\phi|_0{\dt}^{k_\beta}+K_c\sum_{\bar\beta\in \bar
I}|\del_t^{\bar\beta_0}D^{\bar\beta'}\phi|_0{\dx}^{\bar k_{\bar\beta}}.  
\smallskip
\end{align*} 
The corresponding version of (S3) is obtained by plugging $\phi_\e$
into (S3') and using the estimates on its derivatives. The result is
\begin{align*} 
&E_1(\tilde K, h, \e) =E_2(\tilde K, h, \e)\\
& =\tilde K
  K_c\sum_{\beta\in I} \e^{1-2\beta_0-|\beta'|}{\dt}^{k_\beta} +\tilde K
  K_c\sum_{\bar\beta\in \bar I}\e^{1-2\bar\beta_0-|\bar\beta'|}
  {\dx}^{\bar k_{\bar\beta}}.
\end{align*} 
From this formula we see that the dependence in the small parameter
  $\e$ is bad since all the exponents of $\e$ are negative, while the
  dependence on $\dt$, $\dx$ is good since their exponents are
  positive.
\end{remark}

\begin{remark}
Assumption (S1) contains two different kinds of information. First, by
taking $\phi\equiv0$ it implies that the scheme is nondecreasing with
respect to the $[u]$ argument. Second, by taking $u\equiv v$ it indicates
that a parabolic equation -- an equation with a $u_t$ term -- is being
approximated. Both these points play a crucial role in the proof of
the comparison principle for \eqref{S} (Lemma \ref{compsche} below).

To better understand that assumption (S1) implies parabolicity of the
scheme, consider the following more restrictive assumption:

\medskip
\noindent {\bf (S1')} {\bf (Monotonicity)}  There exists $\lambda\geq
0$, $\bar K>0$ such that if $u\leq v$, $u,v\in C_b(\G_h)$, and
$\phi:[0,T]\ra\R$ is smooth, then
\begin{align*}
&S(h,t,x,r+\phi(t),[u+\phi]_{t,x})\\
& \geq  S(h,t,x,r,[v]_{t,x})
+\phi'(t)-\bar K\dt|\phi''|_0-\lambda\phi^+(t) \quad \text{in} \quad \G_h^+.
\medskip
\end{align*}
Here $h=(\dt, h')$ where $h'$ representing a small parameter
related to e.g. the space discretization.
It is easy to see that (S1') implies (S1), e.g. with the same
value for $\lambda$ and the following values of $\mu$ and $h_0$:
$$\mu =\lambda+1 \qquad\text{and}\qquad h_0^{-1} = 2\bar K
e^{(\lambda+1)T}(\lambda+1)(2+(\lambda+1)T).$$   
Assumption (S1') is satisfied for all {\em monotone} finite difference in
time approximations of \eqref{E}, e.g. {\em  monotone} Runge-Kutta methods and
{\em monotone} multi-step methods, both explicit and implicit
methods. We have emphasized the word monotone because whereas many
Runge Kutta methods actually lead to monotone schemes
for \eqref{E} (possibly under a CFL condition), it seems that the most commonly
used multistep methods (Adams-Bashforth, BDS) do not. We refer to
\cite{Shu} for an example of a multistep method that yields a monotone
approximation of \eqref{E}. 
\end{remark}

\medskip
\noindent{\bf Proof of Theorem~\ref{mainres}.}
We start by proving that conditions (S1) and (S2) imply a comparison
result for bounded continuous sub and supersolutions of \eqref{S}.
\begin{lemma}\label{compsche}
Assume (S1), (S2), and that $u,v\in C_b(\G_h)$ satisfy
$$ S(h,t,x,u (t,x),[u]_{t,x})\leq g_1 \quad\hbox{in  } \G_h^+\; ,$$
$$ S(h,t,x,v (t,x),[v]_{t,x})\geq g_2 \quad\hbox{in  } \G_h^+\; ,$$
where $g_1, g_2 \in C_b(\G_h)$. Then
$$u-v\leq e^{\mu t}|(u(0,\cdot)-v(0,\cdot))^+|_0 + 2te^{\mu
  t}|(g_1-g_2)^+|_0,$$ 
where $\lambda$ and $\mu$ are given by (S1).
\end{lemma}

\begin{proof}
1. First, we notice that it suffices to prove the lemma in the case
\begin{align}
\label{AIV}
u(0,x)-v(0,x)\leq 0 \quad\text{in}\quad\G_h^0,
\end{align}
\begin{align}
\label{AST}
g_1 (t,x)-g_2(t,x)\leq 0 \quad\text{in}\quad\G_h.
\end{align}
The general case follows from this result after noting that, by (S1),
$$w=v+e^{\mu t}\left( |(u(0,\cdot)-v(0,\cdot))^+|_0 + 2t
|(g_1-g_2)^+|_0\right)\; ,$$ 
satisfies $S(h,t,x,w (t,x),[w]_{t,x})\geq g_1$ in $\G_h^+$ and
$u(0,x)-w(0,x)\leq0$ in $\G_h^0$.

\noindent 2. We assume that (\ref{AIV}) and (\ref{AST}) hold and, for
$b\geq 0$, we set $\psi_b (t)=e^{\mu t}2b t$ where $\mu$ is given by
(S1) and 
$$ M(b)=\sup_{\G_{h}}\{u-v-\psi_b\}\; .$$

We have to prove that $M(0) \leq 0$ and we argue by contradiction
assuming that $M(0) > 0$. 

\noindent 3. First we consider some $b\geq 0$ for which $M(b) >0$ and
take a sequence $\{(t_n,x_n)\}_n \subset \G_h$ such that
$$\delta_n:=M(b)-(u-v-\psi_b)(t_n,x_n)\ra 0\quad \text{as}\quad
  n\ra\infty.$$ 
Since $M(b)>0$ and \eqref{AIV} holds, $t_n>0$ for all
  sufficiently large $n$ and for such $n$, we have
\begin{align*}
g_1&\geq
S(h,t_n,x_n,u,[u]_{t_n,x_n}) && \text{\small ($u$
  subsolution)}\\ 
&\geq
S(h,t_n,x_n,v+\psi_b +M(b)-\delta_n,[v+\psi_b+M(b)]_{t_n,x_n})
&&\text{\small (S1), $\phi\equiv0$}\\ 
&\geq \omega(\delta_n)\\
&\quad+S(h,t_n,x_n,v+\psi_b+M(b),[v+\psi_b+M(b)]_{t_n,x_n}) &&
\text{\small (S2)}\\ 
&\geq \omega(\delta_n) + b -\lambda M(b) +
S(h,t_n,x_n,v,[v]_{t_n,x_n}) && 
\text{\small (S1), $\phi=\psi+M$}\\
&\geq \omega(\delta_n)+ b -\lambda M(b) +g_2,  &&  \text{\small ($v$
  supersolution)}
\end{align*}
where we have dropped the dependence in $t_n,x_n$ of $u$, $v$ and
  $\psi_b$ for the sake of simplicity of notation.
Recalling \eqref{AST} and sending $n\ra\infty$ lead to
$$ b -\lambda M(b) \leq 0\; .$$

\noindent 4. Since $M(b) \leq M(0)$, the above inequality yields a
contradiction for $b$ large, so for such $b$, $M(b) \leq 0$. On the
other hand, since $M(b)$ is a continuous 
function of $b$ and $M(0) >0$, there exists a minimal solution $\bar b
>0$ of $M(\bar b) =0$. For $\delta >0$ satisfying $\bar b-\delta >0$,
we have $M(\bar b-\delta) >0$ and 
$M(\bar b-\delta)\to 0$ as $\delta \to 0$. But, by 3 we have
$$ \bar b-\delta \leq \lambda M(\bar b-\delta) ,$$
which is a contradiction for $\delta$ small enough since $\bar b >0$.
 \end{proof}

Now we turn to the {\bf proof of the upper bound}, i.e. of (a). We
just sketch it since it relies on the regularization procedure of
Krylov which is used in Section~\ref{Sec:Sw}. We also refer to Krylov
\cite{Kr:HJB1,Kr:HJB2}, Barles and Jakobsen \cite{BJ:Rate,J:Par} for
more details. The main steps are: \smallskip 

\noindent 1. Introduce the solution $u^\eps$ of
\begin{align*}
u^{\eps}_ t+\sup_{0 \leq s \leq\eps^2, |e|\leq\eps}\, F(t+s,x+e,
u^{\eps}(t,x) , Du^{\eps}, D^2 u^{\eps}) &= 0 && \hbox{in
}Q_{T+\eps^2},\\
 u^{\eps}(x,0) &= u_0 (x) && \hbox{in  }\R^N.
\end{align*}
Essentially as a consequence of Proposition \ref{WPi}, it follows that $u^\eps$
belongs to $\CC^{0,1}(Q_T)$ with a uniform $\CC^{0,1}(Q_T)$-bound 
$\bar K$. \smallskip

\noindent 2. By analogous arguments to the ones used in Section~\ref{Sec:Sw},
it is easy to see that $u_{\eps}:=u^{\eps}*\rho_\eps$ is a subsolution
of \eqref{E}. By combining regularity and continuous dependence
results (Theorem \ref{CD} in the Appendix), we also have
$|u_{\eps}-u|_0 \leq C\eps$ where $C$ only depends $T$ and $K$ in
(A1).\smallskip 

\noindent 3. Plugging $u_{\eps}$ into the scheme and using (S3)(i) and the
  uniform estimates on $u^{\eps}$ we get
$$ S(h,t,x,u_{\eps} (t,x),[u_{\eps}]_{t,x})\leq E_1(\bar K,h,\eps)
\quad\hbox{in } \G_h^+\; ,$$ where $\bar K$ is the above mentioned
$\CC^{0,1}$--uniform estimate on $u^{\eps}$ which depends only on the data
and is essentially the same as for $u$. \smallskip

\noindent 4. Use Lemma~\ref{compsche} to compare $u_{\eps}$ and
$u_h$ and conclude by using the control we have on $u-u_{\eps}$ and
  by taking the minimum in $\eps$.

\smallskip

We now provide the {\bf proof of the lower bound}, i.e. of
(b). Unfortunately, contrarily to the proof of (a), we do not know how
to obtain a sequence of approximate, global, smooth supersolutions. As
in \cite{BJ:Rate2}, we are going to obtain approximate ``almost
smooth'' supersolutions which are in fact supersolutions which are
smooth at the "right points". We build them by considering the
following switching system approximation of \eqref{E}:
\begin{align}
\label{Sw3}
F^\eps_i(t,x,v^\eps,\del_t v^\eps_i,Dv^\eps_i,D^2v^\eps_i)&=0 &&\!\!\!
\text{in}\ \ Q_{T+2\eps^2}, \ i\in\II:=\{1,\dots,M\}, \\
v^\eps(0,x)&=v_0(x) &&\!\!\!\text{in}\ \ \R^N,\nonumber
\end{align}
where $v^\eps = (v^\eps_1, \cdots, v^\eps_M)$, $v_0=(u_0,\dots,u_0)$, 
\begin{align}\label{esc}
&F^\eps_i(t,x,r,p_t,p_x,X)=
\\ &\max\Big\{p_t+\min_{0\leq s \leq
\eps^2, |e|\leq\eps}\LL^{\alp_i}(t + s - \eps^2 ,x + e,r_i,p_x,X);\
r_i-\M_ir\Big\},\nonumber
\end{align}
and $\LL$ and $\M$ are defined below \eqref{E} and \eqref{Sw1}
respectively. The solution of this system is expected to be close to
the solution of \eqref{E} if $k$ and $\eps$ are small and
$\{\alp_i\}_{i\in\II}\subset\A$ is a sufficiently refined grid for
$\A$. This is where the assumption (A2) plays a role.

For equation \eqref{Sw3}, we have the following result.
\begin{lemma}
\label{l:sw}
Assume (A1).

(a) There exists a unique solution $v^\eps$ of \eqref{Sw3} satisfying
$|v^\eps|_1\leq \bar K$, where $\bar K$ only depends on $T$ and $K$ from (A1).

(b) Assume in addition (A2) and let $u$ denote the solution of
\eqref{E}. For $i\in\II$, we introduce the functions $\bar v_{i}^\eps : [-\eps^2, T + \eps^2]\times \R^N\ra
\R$  defined by 
$$\bar v_{i}^\eps (t,x):= v_{i}^\eps (t-\eps^2,x).$$ 
Then, for any $\delta>0$, there are $M\in\N$ and
$\{\alp_i\}_{i=1}^M\subset\A$ such that
$${\rm max}_{i}\,|u-\bar v_{i}^\eps|_0\leq
C(\eps+k^{1/3}+\delta),$$
where $C$ only depends on $T$ and $K$ from (A1).
\end{lemma}

In order to simplify the arguments of the proof of the lower bound (to
have the simplest possible formulation of Lemma~\ref{wtsupersol}
below), we need the solutions of the 
equation with ``shaken coefficients'' to be defined in a slightly larger
domain than $Q_T$. More precisely on 
$$Q_T^\eps:=(-\eps^2, T + \eps^2]\times\R^N.$$
This is the role of the $\bar 
v_{i}^\eps$'s. In fact they solve the same system of equations as
the $v_{i}^\eps$'s but on $Q^\eps_T$ and with
$\LL^{\alp_i}(t + s - \eps^2 ,x + e,r_i,p_x,X)$ being replaced by 
$\LL^{\alp_i}(t + s ,x + e,r_i,p_x,X)$ in \eqref{esc}. 

The (almost) smooth supersolutions of \eqref{E} we are looking for are
built out of the $\bar v_{i}^\eps$'s by mollification. Before giving
the next lemma, we remind the reader that the sequence of mollifiers
$\{\rho_\eps\}_{\eps}$ is defined at the end of the introduction.

\begin{lemma}
\label{l:w}
Assume (A1) and define ${v_\eps}_i:=\rho_\eps* \bar v^\eps_i : \ol
Q_{T+\eps^2}\ra 
\R$ for $i\in\II$. 

(a) There is a constant $C$ depending only on 
$T$ and $K$ from (A1), such that 
$$|{v_\eps}_{j}-\bar  v^\eps_i|\leq C(k+\eps)\quad \text{in} \quad 
Q_{T+\eps^2},\quad i,j\in\II.$$

(b) Assume in addition that $\eps \leq (8\sup_i[v^\eps_i]_{1})^{-1}k$. For
every $(t,x)\in Q_T$, if $j:=\argmin_{i\in\II}{v_\eps}_i(t,x)$, then 
\begin{equation*}
  \del_t{v_\eps}_j(t,x) +
    \LL^{\alp_j}(t,x,{v_\eps}_j(t,x),D{v_\eps}_j(t,x), D^2 {v_\eps}_j(t,x)) 
    \geq 0.
\end{equation*}
\end{lemma}

The proofs of these two lemmas will be given at the end of this section.

The key consequence is the following result which is the corner-stone
of the proof of the lower bound.
\begin{lemma}
\label{wtsupersol} 
Assume (A1) and that $\eps \leq (8\sup_i[v^\eps_i]_{1})^{-1}k$. Then the
function $w:=\min_{i\in\II}{v_\eps}_i$ is an approximate supersolution
of the scheme \eqref{S} in the sense that
$$S(h,t,x,w (t,x),[w]_{t,x})\geq - E_2 (\bar K, h, \eps)
\quad\hbox{in }\G_h^+,$$ 
where $\bar K$ comes from Lemma~\ref{l:sw}.
\end{lemma}

\begin{proof} Let $(t,x)\in Q_T$ and $j$ be as in Lemma~\ref{l:w}
  (b). We see that $w(t,x) = {v_\eps}_j(t,x)$ and $w \leq {v_\eps}_j$ in
  $\G_h$, and hence the monotonicity of the scheme (cf. (S1)) implies that
$$ S(h,t,x,w (t,x),[w]_{t,x}) \geq  S(h,t,x,{v_\eps}_j
  (t,x),[{v_\eps}_j]_{t,x}).$$
But then, by (S3)(ii),
\begin{align*}
&S(h,t,x,w (t,x),[w]_{t,x})  \\
&\geq   \del_t{v_\eps}_j(t,x) + \LL^{\alp_j}(t,x,{v_\eps}_j(t,x),
D{v_\eps}_j(t,x), D^2 {v_\eps}_j(t,x)) - E_2 (\bar K, h, \eps),
\end{align*}
and the proof complete by applying Lemma~\ref{l:w} (b). 
\end{proof}

It is now straightforward to conclude the proof of the lower bound, we
simply choose $k= 8\sup_i[v^\eps_i]_{1} \eps$ and use
Lemma~\ref{compsche} to compare $u_h$ and $w$. This yields
$$ u_h-w  \leq e^{\mu t}|(u_{h,0}-w(0,\cdot))^+|_0 + 2te^{\mu
  t}E_2 (\bar K, h, \eps) \quad\hbox{in  }\G_h .$$ 
But, by Lemmas \ref{l:sw} (b) and \ref{l:w} (a), we have
$$ |w - u|_0 \leq C(\eps+k+k^{1/3}+\delta),$$
and therefore 
$$ u_h- u  \leq e^{\mu t}|(u_{h,0}-u_0)^+|_0 + 2te^{\mu
  t}E_2 (\bar K, h, \eps) + C(\eps+k+k^{1/3}+\delta) \quad\hbox{in  }\G_h ,$$ 
for some constant $C$. In view of our choice of $k$, we conclude the proof by
minimizing w.r.t $\eps$.

\bigskip
Now we give the proofs of Lemmas \ref{l:sw} and \ref{l:w}.

\begin{proof}[Proof of Lemma \ref{l:sw}]
{\bf 1.} We first approximate \eqref{E} by
\begin{align*}
v_t+\sup_{i\in\II}\LL^{\alp_i}(t,x,v,Dv,D^2v)&=0&&\text{in}\quad Q_T,\\
v(0,x)&=u_0(x)&&\text{in}\quad \R^N.
\end{align*}
From assumption (A2) and Lemmas \ref{WP} and \ref{CD} in the Appendix,
it follows that there exists a unique solution $v$ of the
above equation satisfying
$$|v-u|_0\leq C\delta,$$
where $C$ only depends on $T$ and $K$ from (A1).

\noindent{\bf 2.} We continue by approximating the above equation by
the following switching system
\begin{align*}
\max\Big\{\del_tv_i+\LL^{\alp_i}(t,x,v_i,Dv_i,D^2v_i);\
v_i-\M_iv\Big\}&=0 && \text{in}\quad Q_T,\quad  i\in\II,\\
v(0,x)&=v_0(x)&&\text{in}\quad\R^N,
\end{align*}
where $v_0=(u_0,\dots,u_0)$ and $\M$ is defined below
\eqref{Sw1}. From Proposition \ref{WPi} 
and Theorem \ref{sw-rate} in Section~\ref{Sec:Sw} we have existence and uniqueness of a solution $\bv=(\bv_1, \cdots, \bv_M)$ of the above system satisfying
$$|\bv_i-v|_0\leq Ck^{1/3},\quad i\in\II,$$
where $C$ only depends on the mollifier $\rho$, $T$, and $K$ from (A1).

\noindent{\bf 3.} The switching system defined in the previous step is
nothing but \eqref{Sw3} with $\eps=0$ or \eqref{Sw2} with the
$\A_{i}$'s being singletons.  Theorems \ref{WP} and \ref{CD} in the
Appendix yield 
the existence and uniqueness of a solution $v^\eps:\ol Q_{T+2\eps^2}\ra\R$ of
\eqref{Sw3} satisfying
$$|v^\eps|_1+\frac{1}{\eps}|v^\eps-\bv|_0\leq C,$$
where $C$ only depends on $T$ and $K$ from (A1).

\noindent{\bf 4.} The proof is complete by combining the estimates in 
steps 1 -- 3, and noting that $|v^\eps-\bar v^\eps|\leq
[v^\eps]_1\eps$ in $Q_{T+\eps^2}$ and (A2) is only needed in step 1. 
\end{proof}

\begin{proof}[Proof of Lemma \ref{l:w}]
We start by (a). From the properties of mollifiers and the H{\"o}lder  
continuity of $\bar v^\eps$, it is immediate
that
\begin{align}
\label{pf1}
|{v_\eps}_i-\bar v^\eps_i|\leq C\eps \quad\text{in}\quad \ol
 Q_{T+\eps^2},\quad i\in\II,
\end{align}
where $C=2\max_i[\bar v^\eps_i]_1= 2\max_i[v^\eps_i]_1$ depends only
on $T$ and $K$ from (A1). Furthermore as we pointed out after the
statement of Lemma~\ref{l:sw}, $\bar v^\eps$ solves a switching system
in $Q^\eps_T$, so  arguing as in the proof of Theorem \ref{sw-rate} in
Section~\ref{Sec:Sw} leads to
$$0\leq \max_i \bar v^\eps_i -\min_i \bar v^\eps_{i}\leq k \quad
\text{in}\quad Q_T^\eps.$$
 From these two estimates, (a) follows.

Now consider (b). Fix an arbitrary point $(t,x)\in
Q_T$ and set
$$j=\argmin_{i\in\II} \, {v_{\eps i}} (t,x).$$
Then, by definition of $\M$ and $j$, we have
$${v_{\eps j}} (t,x)-\M_j{v_\eps}(t,x)= \max_{i\neq 
j}\left\{{v_{\eps j}}(t,x) - {v_{\eps i}}(t,x)
-k\right\}\leq -k, $$ 
and the bound \eqref{pf1} leads to
$$\bar v^\eps_j(t,x)-\M_j \bar v^\eps(t,x)\leq -k+2\max_i[v^\eps_i]_12\eps .$$ 
Next, by using the H{\"o}lder continuity of $\bar v^\eps$ (Lemma
\ref{l:sw}), for any $(\tb,\xb)\in Q_T^\eps$, we have 
$$\bar v^\eps_j (\tb,\xb)-\M_j \bar v^\eps(\tb,\xb) \leq
 -k+2\max_i[v^\eps_i]_1(2\eps+|x-\xb|+|t-\tb|^{1/2}).$$
 From this we conclude that if $|x-\xb|<\eps$, $|t-\tb|<\eps^2$, and
 $\eps \leq (8\max_i[v^\eps_i]_1)^{-1}k$, then
$$\bar v^\eps_j (\tb,\xb) -\M_j \bar v^\eps (\tb,\xb) <0,$$
and by equation \eqref{Sw3} and the definition of $\bar v^\eps$, $\bar
v^\eps(t,x)=v^\eps(t-\eps^2,x)$, 
$$\del_t \bar v^\eps_j (\tb,\xb) +\inf_{0\leq s \leq \eps^2,
|e|\leq\eps}\LL^{\alp_j}(\tb + s,\xb + e, \bar  v^\eps_j(\tb,\xb),
D \bar v^\eps_j(\tb,\xb), D^2 \bar v^\eps_j(\tb,\xb))=0.$$ 
After a change of variables, we see that, for every $0\leq s <
\eps^2, |e|<\eps$, 
\begin{align}
\label{ineq5}
&\del_t \bar  v^\eps_j(t - s,x - e) (t,x)\\
&+\LL^{\alp_j}(t,x, \bar v^\eps_j(t - s,x -
e),D \bar  v^\eps_j(t - s,x - e),D^2 \bar  v^\eps_j(t - s,x - e))\geq 0.\nonumber
\end{align}
In other words, for every
$0\leq s < \eps^2, |e|< \eps$, $\bar v^\eps_j(t - s,x - e)$ is a (viscosity)
supersolution at $(t,x)$ of  
\begin{align}
\label{lin}
\chi_t+\LL^{\alp_j}(t,x,\chi,D\chi,D^2\chi)=0.
\end{align}

By mollifying \eqref{ineq5} (w.r.t. the $(s,e)$-argument) we see that
${v_\eps}_j$ is also a smooth supersolution of \eqref{lin} at $(t,x)$
and hence a  
(viscosity) supersolution of the HJB equation \eqref{E} at $(t,x)$. 
This is correct since ${v_\eps}_j$ can be viewed as the limit of convex
combinations of supersolutions $\bar v^\eps_j(t - s,x - e)$ of the linear and hence
concave equation \eqref{lin}, we refer to the proof of Theorem
\ref{sw-rate} and to the Appendix in \cite{BJ:Rate} for the details. 
We conclude the proof by noting that since ${v_\eps}_j$ is smooth, it
is in fact a classical supersolution of \eqref{E} at $x$. 
\end{proof}

\section{Monotone Finite Difference Methods}
\label{Sec:FDM}

In this section, we apply our main result to finite difference
approximations of \eqref{E} based on the $\vartheta$-method
approximation in time and two different approximations in space: One
proposed by Kushner \cite{KD:Book} which is monotone when 
$a$ is diagonal dominant and a (more) general approach based on
directional second derivatives proposed by Bonnans and Zidani
\cite{BZ}, but see also Dong and Krylov \cite{DK:ConstCoeff}.
For simplicity we take $h=(\dt,\dx)$ and consider the uniform grid
$$\G_h=\dt\{0,1,\dots,n_T\}\times\dx \Z^N.$$

\subsection{Discretization in space}

To explain the methods we first write equation \eqref{E} like 
\begin{align*}
u_t+\sup_{\alpha\in\A}\Big\{-L^\alpha u - c^\alpha(t,x) u -
f^\alpha (t,x)\Big\}=0 \quad\text{in}\quad Q_T,
\end{align*}
where
\begin{align*}
L^\alpha \phi(t,x)= \tr [a^\alpha(t,x) D^2\phi(t,x)] + b^\alpha(t,x) D
\phi(t,x).
\end{align*}
To obtain a discretization in space we approximate $L$ by a finite difference
operator $L_h$, which we will take to be of the form \begin{align}
\label{FDM}
L_h^\alpha\phi(t,x) = \sum_{\beta \in\CS}C^\alp_h(t,x,
\beta)(\phi(t,x+\beta \dx)-\phi(t,x)), 
\end{align}
for $(t,x)\in\G_h$, where the {\em stencil} $\CS$ is a finite subset of
$\Z^N\setminus\{0\}$, and where
\begin{align}
\label{pos}
C^\alp_h(t,x,\beta)\geq 0\quad\text{for all } \beta \in\CS, (t,x)\in
  \G_h^+, h=(\dx,\dt)>0, \alpha\in\A. 
\end{align}
The last assumption says that the difference approximation is of {\em
  positive} type. This is a sufficient assumption for monotonicity in
  the stationary case.\\

(i) \underline{The approximation of Kushner.}\\ We denote by
$\{e_i\}_{i=1}^N$ the standard basis in $\R^N$ and define
\begin{align}
\label{KuAppr}
L_h^\alpha\phi = &\sum_{i=1}^N
\Big[\frac{a^\alpha_{ii}}{2} \Delta_{ii} + \sum_{j\neq i}
\Big( \frac{a^{\alpha+}_{ij}}{2} \Delta_{ij}^+-
\frac{a^{\alpha-}_{ij}}{2} \Delta_{ij}^- \Big)+ b_i^{\alpha+} \delta^+_i
- b_i^{\alpha-}\delta^-_i \Big]\phi, 
\end{align}
where $b^+=\max\{b,0\}$, $b^-=(-b)^+$ ($b=b^+-b^-$), and
\begin{align*}
\delta^{\pm}_{i} w(x)&=\pm\frac{1}{\dx}\{w(x \pm e_i\dx)-w(x) \},\\
\Delta_{ii} w(x)&= \frac{1}{\dx^2}\{w(x+e_i \dx) -2w(x) + w(x- e_i
  \dx)\},\\ 
\Delta^+_{ij} w(x)&= \frac{1}{2\dx^2}\{2w(x) + w(x+e_i\dx+e_j\dx)
+ w(x-e_i\dx-e_j\dx)\} \\ 
&\hspace{-1cm} - \frac{1}{2\dx^2}\{w(x+e_i\dx) + w(x-e_i\dx) + w(x+e_j\dx) +
w(x-e_j\dx)\},\\
\Delta^-_{ij} w(x) & = -\frac{1}{2\dx^2}\{2w(x) + w(x+e_i\dx-e_j\dx) +
w(x-e_i\dx+e_j\dx)\}\\ 
&\hspace{-1cm} +\frac{1}{2\dx^2}\{w(x+e_i\dx) + w(x-e_i\dx) +
w(x+e_j\dx) + w(x-e_j\dx)\} .
\end{align*}
The stencil is $\CS=\{\pm e_i, \pm(e_i\pm e_j):i,j=1,\dots, N\}$, and it
is easy to see that the coefficients in \eqref{FDM} are
\begin{align*}
&C^\alp_h(t,x,\pm e_i)\hspace{-1.7cm}&&= \frac{a^\alpha_{ii}(x)}{2\dx^2} -
\sum_{j\neq i} \frac{|a^\alpha_{ij}(x)|}{4\dx^2} +
\frac{b^{\alpha\pm}_i(x)}{\dx},\\ 
&C^\alp_h(t,x,e_ih \pm e_jh)\hspace{-1.7cm}&&=
\frac{a^{\alpha\pm}_{ij}(x)}{2\dx^2}, \quad i\neq j,\\
&C^\alp_h(t,x,-e_ih \pm e_jh)\hspace{-1.7cm}&&=
\frac{a^{\alpha\mp}_{ij}(x)}{2\dx^2},\quad i\neq j. 
\end{align*}
The approximation is of positive type \eqref{pos} if and only if $a$
is diagonal dominant, i.e. 
\begin{align}
\label{DiagDom}
a^\alpha_{ii}(t,x) - \sum_{j\neq i} |a^\alpha_{ij}(t,x)| \geq 0 \quad
\text{in}\quad Q_T, \quad \alp\in\A, \quad i=1,\dots,N.
\end{align}

(ii) \underline{The approximation of Bonnans and Zidani.}\\ We assume
that there is a (finite) stencil $\bar\CS\subset\Z^N\setminus\{0\}$
and a set of positive coefficients $\{\bar a_\beta : \beta
\in\bar\CS\}\subset \R_+$ such that
\begin{align}
\label{GenDom}
a^\alpha(t,x)=\sum_{\beta\in\bar\CS}\bar a^\alpha_\beta (t,x)\beta^T\beta
 \quad  \text{in}\quad Q_T,\quad \alpha \in \A.
\end{align}
Under assumption \eqref{GenDom} we may rewrite the operator $L$ using
second order directional derivatives $D^2_\beta=\tr[\beta\beta^TD^2]=(\beta
\cdot D)^2$,
$$L^\alpha \phi(t,x)= \sum_{\beta\in\bar\CS}\bar a^\alpha_\beta(t,x)
D^2_{\beta}\phi(t,x) + b^\alpha(t,x) D \phi(t,x).$$
The approximation of Bonnans and Zidani is given by
\begin{align}
\label{BoAppr}
L_h^\alpha\phi = &\sum_{\beta\in\bar\CS}\bar a^\alpha_\beta
\Delta_\beta\phi +\sum_{i=1}^N\Big[ b_i^{\alpha+} \delta^+_i
- b_i^{\alpha-}\delta^-_i \Big]\phi, 
\end{align}
where $\Delta_{\beta}$ is an approximation of $D^2_\beta$ given by
\begin{align*}
&\Delta_{\beta} w(x)= \frac{1}{|\beta|^2\dx^2}\{w(x+\beta \dx) -2w(x)
  + w(x- \beta \dx)\}.
\end{align*}
In this case, the stencil is $\CS=\pm\bar\CS\cup\{\pm e_i:
i=1,\dots,N\}$ and the coefficients corresponding to \eqref{FDM} are given by 
\begin{align*}
C^\alp_h(t,x,\pm e_i)&= \frac{b^{\alpha\pm}_i(x)}{\dx},&& i=1,\dots,N,\\
C^\alp_h(t,x,\pm\beta)&= \frac{\bar
  a^\alpha_\beta (t,x)}{|\beta|^2\dx^2},&&\beta \in\bar\CS,
\end{align*}
and the sum of the two whenever $\beta=e_i$. 
Under assumption \eqref{GenDom}, which is more general than
\eqref{DiagDom} (see below), this approximation is always of
positive type.\\

For both approximations there is a constant $C>0$, independent of $\dx$, such
that, for every $\phi\in C^4(\R^N)$ and $(t,x)\in\G^+_h$, 
\begin{align}
\label{Consist}
|L^\alpha\phi(t,x)-L_h^\alpha\phi(t,x)|\leq C(|b^\alpha|_0|D^2\phi|_0\dx +
 |a^\alpha|_0|D^4\phi|_0\dx^2).
\end{align}

\subsection{The fully discrete scheme} To obtain a fully discrete scheme,
we apply the $\vartheta$-method, $\vartheta\in[0,1]$, to discretize
the time derivative. The result is the following scheme,
\begin{align}
\label{FFDM}
u(t,x)=&\ u(t-\dt,x)\\
&-(1-\vartheta)\dt\sup_{\alpha\in\A}\{-L^\alpha_h
u - c^\alpha u - 
f^\alpha \}(t-\dt,x)\nonumber\\
&-\vartheta\dt\sup_{\alpha\in\A}\{-L^\alpha_h
u - c^\alpha u - f^\alpha \}(t,x)
 \quad\text{in}\quad \G^+_h.\nonumber
\end{align}
The case $\vartheta=0$ and $\vartheta=1$ correspond to the forward and
backward Euler time-discretizations respectively, while for
$\vartheta=1/2$ the scheme is a generalization of the second order in
time Crank-Nicholson scheme. Note that the scheme is implicit except
for the value $\vartheta=0$. We may write \eqref{FFDM} in the form
\eqref{S} by setting
\begin{align*}
&S(h,t,x,r,[u]_{t,x})=\sup_{\alpha\in\A}\Big\{\Big[\frac1\dt-\vartheta
c^\alpha+\vartheta\sum_{\beta\in\CS}C^\alp_h(t,x,\beta)\Big]r\\
&-\Big[\frac1\dt+(1-\vartheta)
c^\alpha -
(1-\vartheta)\sum_{\beta\in\CS}C^\alp_h(t,x,\beta)\Big][u]_{t,x}(-\dt,0)\\  
&-\sum_{\beta\in\CS}C^\alp_h(t,x,\beta)\Big[\vartheta[u]_{t,x}(0,\beta\dx) +
  (1-\vartheta)[u]_{t,x}(-\dt,\beta\dx)\Big]\Big\},
\end{align*}
where $[u]_{t,x}(s,y)=u(t+s,x+y)$. Under assumption \eqref{pos} 
the scheme \eqref{FFDM} is monotone (i.e. satisfies
(S1) or (S1')) provided the following CFL conditions hold
\begin{align}
\label{CFL}
\dt\,(1-\vartheta)
\Big(-c^{\alpha}(t,x)+\sum_{\beta\in\CS} 
C^\alp_h(t,x,\beta)\Big)&\leq 1,\\
\label{CFL2}
\dt \,\vartheta
\Big(c^\alpha(t,x)-\sum_{\beta\in\CS}C^\alp_h(t,x,\beta)\Big)&\leq 1. 
\end{align}
Furthermore, in view of (A1) and \eqref{Consist}, Taylor expansion in
\eqref{FFDM} yields the following consistency result for smooth functions 
$\phi$ and $(t,x)\in\G_h^+$,
\begin{align*}
&|\phi_t+F(t,x,\phi,D\phi,D^2\phi) - S(h,t,x,\phi,[\phi]_{t,x})|\\
&\leq
C(\dt|\phi_{tt}|_0+\dx|D^2\phi|_0+\dx^2|D^4\phi|_0+(1-\vartheta)\dt(|D\phi_t|_0+|D^2\phi_t|_0)).\nonumber 
\end{align*}
The $(1-\vartheta)\dt$-term is a non-standard term coming from the
fact that we need the equation and the scheme to be satisfied in the
{\em same} point, see assumption (S3). The necessity of this
assumption follows from the proof of Theorem~\ref{mainres}. 

We have seen that if \eqref{pos} and \eqref{Consist} hold along with
the CFL conditions \eqref{CFL} and \eqref{CFL2} then the scheme
\eqref{FFDM} satisfies assumptions (S1) -- (S3) in
Section~\ref{Sec:HJB}. Theorem 
\ref{mainres} therefore yields the following error bound:
\begin{theorem}\label{rate1}
Assume (A1), (A2), \eqref{pos}, \eqref{Consist}, \eqref{CFL},
\eqref{CFL2} hold. If $u_h\in C_b(\G_h)$ is the solution of
\eqref{FFDM} and $u$ is the solution of \eqref{E}, then there is $C>0$
such that in $\G_h$ 
\begin{align*}
-e^{\mu t}|(u_0-u_{0,h})^-|_0-C|h|^{\frac1{5}}
\leq u-u_h
\leq e^{\mu t}|(u_0-u_{0,h})^+|_0+C|h|^{\frac1{2}},
\end{align*}
where $|h| := \sqrt{\dx^2+\dt}$.
\end{theorem}

\begin{remark}
Except when $\vartheta=1$, the CFL condition \eqref{CFL}
essentially implies that $\dt \leq C\dx^2$. Therefore $\dt$ and $\dx^2$ play
essentially the same role. Also note that the CFL condition
\eqref{CFL2} is satisfied if e.g. $\dt\leq
(\sup_\alpha|(c^\alpha)^+|_0)^{-1}$.  
\end{remark}
\begin{remark}
Even though the above consistency relationship is not quite the ``standard''
one, it gives the correct asymptotic behavior of our scheme. First
of all note that the new term, the $(1-\vartheta)$-term, behaves just
like the $\Delta t$ and  ${\Delta x}^2$ terms. To see this, we note
that according to (S3) we 
only need the above relation when $\phi$ is replaced by $\phi_\eps$
defined in (S3). But for $\phi_\eps$ we have $|\phi_{\eps,
  tt}|_0\approx|D^4\phi_{\eps}|_0\approx|D^2\phi_{\eps, t}|_0\approx
\tilde K \eps^{-3}.$ By the CFL conditions \eqref{CFL} and \eqref{CFL2} we have
essentially that $\dx^2\approx\dt$, so  
$$\dt|\phi_{\eps, tt}|_0\approx\dx^2|D^4\phi_{\eps}|_0\approx\dt|D^2\phi_{\eps,
  t}|_0\approx \tilde K \dx^2\eps^{-3}.$$

Next note that for $\vartheta=1/2$ (the Cranck-Nicholson case) the
scheme is formally second order in time. However this is no longer the
case for the monotone version. It is only first order in time due to the CFL 
condition which implies that $\dx^2\|D^4\phi\|=C\dt\|D^4\phi\|$.
\end{remark}

\begin{proof} In this case
\begin{align*}
&E_1 (\bar K,h,\eps) = E_2 (\bar K,h,\eps)\\
&=C(\dt\eps^{-3}+\dx\eps^{-1} + \dx^2\eps^{-3} +
  (1-\vartheta)\dt(\eps^{-2}+\eps^{-3})). 
\end{align*}
So we have to minimize w.r.t. $\eps$ the following functions
$$ \eps + C(\dt\eps^{-3}+\dx\eps^{-1}+\dx^2\eps^{-3}),$$
$$ \eps^{1/3} + C(\dt\eps^{-3}+\dx\eps^{-1}+\dx^2\eps^{-3}).$$
By minimizing separately in $\dt$ and $\dx$, one finds that $\eps$ has
to be like $\dt^{1/4}$ and $\dx^{1/2}$ in the first case, and that
$\eps^{1/3}$ has to be like $\dt^{1/10}$ and $\dx^{1/5}$ in the second
case. The result now follows by taking 
$\eps=\max(\dt^{1/4},\dx^{1/2})$ in the first case and
$\eps^{1/3}=\max(\dt^{1/10},\dx^{1/5})$ in the second case.
\end{proof}

\subsection{Remarks}
For approximations of nonlinear equations monotonicity is a
key property since it ensures (along with consistency) that the
approximate solutions converge to the {\em correct} generalized solution
of the problem (the viscosity solution in our case). This is not the case for
nonmonotone methods, at least not in any generality. 

However, the monotonicity requirement poses certain problems. Monotone
schemes are low order schemes, and maybe more importantly, it is not always
possible to find consistent monotone approximations for a given problem.
To see the last point we note that in general the second derivative
coefficient matrix $a$ is only positive semidefinite, while the monotone
schemes of Kushner and Bonnans/Zidani require the stronger assumptions
\eqref{DiagDom} and \eqref{GenDom} respectively. In fact, in Dong and
Krylov \cite{DK:ConstCoeff} it was proved that if an operator $L$
admits an approximation $L_h$ of the form \eqref{FDM} which is of
positive type, then $a$ has to satisfy \eqref{GenDom} (at least if $a$
is bounded). 

This is a problem in real applications, e.g. in finance, and it was
this problem was the motivation behind the approximation of 
Bonnans and Zidani. First of all we note that their condition
\eqref{GenDom} is more general than \eqref{DiagDom} because any
symmetric $N\times N$ matrix $a$ can be decomposed as  
$$a=\sum_{i=1}^N\sum_{j\neq i} (a_{ii}-|a_{ij}|)e_ie_i^T +
\frac{a_{ij}^+}{2}(e_i+e_j)(e_i+e_j)^T +
\frac{a_{ij}^-}{2}(e_i-e_j)(e_i-e_j)^T,
$$  
where the coefficients are nonnegative if and only if $a$ is
diagonal dominant. More importantly, it turns out that any symmetric positive
semidefinite matrix can be approximated by a sequence of matrices
satisfying \eqref{GenDom}. In Bonnans, Ottenwaelter, and Zidani
\cite{BOZ}, this was proved in the case of symmetric $2\times2$
matrices along with an explicit error bound and an algorithm for
computing the approximate matrices. Because of continuous dependence
results for the equations, convergence of the coefficients immediately
imply convergence of the solutions of the corresponding
equations. Hence the Bonnans/Zidani approximation yields a way of
approximating general problems where $a$ is only positive semidefinite.

\section{Semigroup Approximations and Splitting Methods}
\label{Sec:Appl}

In this section, we consider various approximations of semigroups
obtained by a semi-discretization in time. In order to simplify the
presentation we start by specializing Theorem \ref{mainres} to the
semigroup setting. To be precise we consider one-step in time
approximations of \eqref{E} given by
\begin{align}
\label{Splitt}
u_h(t_n,x)&=S_h(t_n,t_{n-1})u_h(t_{n-1},x)&& \quad{in} \quad \R^N,\\
u_h(0,x)&=u_{h,0}(x)&& \quad{in} \quad \R^N,\nonumber
\end{align}
where $t_0=0 < t_1 <\cdots < t_n < \cdots<t_{n_T}=T$, $h:= \max_n
(t_{n+1}-t_n)$, and the approximation semigroup $S_h$ satisfies the
following sub and superconsistency requirements: There 
exist a constant $K_c$, a subset $I$ of $\N \times
\N^N$, and constants $\gamma_\beta,\delta_\beta$ for $\beta\in I$ such
that for any smooth functions $\phi$, 
\begin{align}
 \label{sconsist}
&\frac1\dt\Big[S_{h}(t_n,t_{n-1})-1\Big]\phi(t_{n-1},x)
-F(t,x,\phi,D\phi,D^2\phi)_{t=t_n}\\
&\leq  K_c\sum_{\beta\in
  I}|\del_t^{\beta_0}D^{\beta'}\phi|^{\gamma_\beta}_0{\dt}^{\delta_\beta},
 \nonumber \\  
\intertext{where $\beta=(\beta_0,\beta')\in I$ for $\beta_0 \in \N$
and $\beta' \in \N^N$, and in a similar way}
&\frac1\dt\Big[S_{h}(t_n,t_{n-1})-1\Big]\phi(t_{n-1},x)
-F(t,x,\phi,D\phi,D^2\phi)_{t=t_n} \label{sconsist2} \\
&\geq  -\bar K_c\sum_{\beta\in
  \bar I}|\del_t^{\beta_0}D^{\beta'}\phi|^{\bar\gamma_\beta}_0
{\dt}^{\bar\delta_\beta}, \nonumber  
\end{align} 
with corresponding data $\bar K_c,\bar I,\bar\gamma_\beta,\bar\delta_\beta$.
We say that the semigroup is {\em monotone} if 
$$\phi\leq\psi \quad \Rightarrow \quad S_{h}(t_n,t_{n-1})\phi\leq
S_{h}(t_n,t_{n-1})\psi, \quad n=1,\dots, n_T,$$
for all continuous bounded functions $\phi,\psi$ for which
$S_h(t)\phi$ and $S_h(t)\psi$ are well defined. 

We have the following corollary to 
Theorem \ref{mainres}. 

\begin{proposition}
\label{semigr_main}
Assume (A1), (A2), and that $S_h$ is a monotone semigroup satisfying
\eqref{sconsist} and \eqref{sconsist2} and which is defined on a subset of
$C_b(\R^N)$. If $u$ is the solution of
\eqref{E} and $u_h$ is the solution of \eqref{Splitt}, then  
$$-C(|u_0-u_{h,0}|_0+{\dt}^{\frac1{10}\wedge r_1})
\leq u-u_h
\leq 
C(|u_0-u_{h,0}|_0+{\dt}^{\frac14\wedge r_2})$$
in $\R^N$, where 
\begin{align*}
r_1 &:= \min_{\beta\in I}\left\{\frac{\delta_\beta}{3(2\beta_0 +
  |\beta'|-1)\gamma_\beta +1}\right\},\\
r_2 &:= \min_{\beta\in \bar I}\left\{\frac{\bar\delta_\beta}{(2\beta_0 +
  |\beta'|-1)\bar\gamma_\beta+1}\right\} ,
\end{align*}
where $|\beta'|$ denotes the sum of the components of $\beta'$.
\end{proposition}

\begin{proof} 
We define
$$S(h,t_n,x,u_h,[u_h]_{t_n,x}) =
\frac1{\dt}\Big(u_h(t_n,x)-[u_h]_{t_{n},x}\Big),$$
where
$$[u_h]_{t_{n},x} = S_h(t_{n-1},t_n)u_h(t_{n-1},x).$$ 
To apply Theorem~\ref{mainres}, we just have to check that (S1) --
(S3) hold and this is clear for (S1) and (S2) (see Remark
\ref{notC_b}). For (S3)(i), note that by \eqref{sconsist} we have
\begin{align*}
&\phi_t+F(t_n,x,\phi,D\phi,D^2\phi)-S(h,t_n,x,\phi,[\phi]_{t_n,x})\\
&\leq \frac12 |\del_t^2\phi|_0\dt+K_c\sum_{\beta\in
  I}|\del_t^{\beta_0}D^{\beta'}\phi|_0^{\gamma_\beta}{\dt}^{\delta_\beta},
\end{align*}
which leads to
$$E_1(\bar K, h, \e)= \frac12 \bar K \eps^{1-4} \dt+
K_c\sum_{\beta\in I} (\bar K
\e^{1-2\beta_0-|\beta'|})^{\gamma_\beta}{\dt}^{\delta_\beta}.$$
The upper bound now follows by optimizing with
respect to $\eps$ as in the proof of Theorem~\ref{rate1}. In a similar
way we may use \eqref{sconsist2} to define $E_2$ and then conclude the
lower bound.
\end{proof}

\begin{remark} In view of the consistency requirements
  \eqref{sconsist} and \eqref{sconsist2}, for schemes like
  \eqref{Splitt} it is natural to think that only the $x$-variable is really
  playing a role and that one can get results on the rate of
  convergence by using this special ``semi-group type''
  structure. More specifically, one might think that a different
  proof using a mollification of the solution with respect to the
  space variable only, can produce the estimates in an easier and
  maybe better way. We tried this strategy but we could not avoid 
  using the short time expansion of the solution of the HJB Equation
  associated with smooth initial data (the short time expansion of
  the semi-group), and this leads to worse rates, even in cases where
  $F$ is smooth. One way of understanding this -- without justifying it
  completely -- consists of looking at our estimates for the
  $\phi_{tt}$-term (cf. (S3)(i) and (ii)). The present approach leads
  to an estimate of 
  order $\eps^{-3}$, while if we use the short time expansion, we are
  lead to a worse estimate of order $\eps^{-4}$. We refer the reader to 
  Subsection~\ref{sds} and in particular to
  Lemma~\ref{ssLem3} below, where short time expansions for semi-groups are
  obtained and used to study the rate of convergence for splitting problems. 
\end{remark}

\subsection{Semidiscrete splitting}\label{sds}

We consider an equation of the form
\begin{align}
\label{spliteq00}
u_t+F_1(D^2u)+F_2(D^2u)=0 \quad\text{in}\quad Q_T,
\end{align}
where
$$F_j(X)=\sup_{\alpha\in\A}\{-\tr[a_j^\alpha X] - f_j^\alpha \},\quad
j=1,2,$$ 
and $a_j^\alpha\geq 0$ are matrices and $f_j^\alpha$ real numbers. We
assume that they are both uniformly bounded in $\alpha$ and are
independent of $(t,x)$. It follows that $F_1$ and $F_2$
are Lipschitz continuous and that (A1) is satisfied.

Let $S$ denote the semigroup of \eqref{spliteq00}, i.e. $S(\dt)\phi$
is the solution at time $t=\dt$ of \eqref{spliteq00} with initial value
$\phi$. Similarly, let $S_1$ and $S_2$ denote the semigroups
associated with the equations $u_t+F_1(D^2u)=0$ and $u_t+F_1(D^2u)=0$. 

We can define a semidiscrete splitting method by taking \eqref{Splitt}
with $t_n:=n\dt$ and 
\begin{align}
\label{SSplitt}
S_h(t_{n-1},t_n)=S_1(\dt)S_2(\dt).
\end{align}
Under the current assumptions all these semigroups map
$W^{1,\infty}(\R^N)$ into itself, they are monotone, and they are
nonexpansive,
$$|\bar S(t)\phi|_0\leq |\phi|_0,$$
for $\phi\in W^{1,\infty}(\R^N)$ and where $\bar S$ denotes one of the
semigroups above.

As soon as we know the consistency relation for this scheme, we can
find an error bound using Theorem \ref{semigr_main}. However,
contrarily to the case of finite different schemes in the previous
section, here the precise form of the consistency requirement is not
well known. We are going to provide such results under different assumptions on
$F_1$, $F_2$.  Our first result is the following:

\begin{lemma} 
\label{ssLem1}
Under the above assumptions, if in addition $|DF_1|\in
W^{1,\infty}(\BS^N)$ and $|DF_2|\in W^{3,\infty}(\BS^N)$, then
\begin{align*}
& - C(\dt |D^2\phi_t|_0 + \dt^2 |D^3\phi|_0^4)- \text{h.o.t.}\\
&\leq \frac1{\dt}[S_h(t)-1]\phi(t_{n-1},x)+F_1(D^2\phi(t_n,x)) +
 F_2(D^2\phi(t_n,x))\\ 
&\leq C(\dt |D^2\phi_t|_0 + \dt |D^3\phi|_0^2)
+ \text{h.o.t.}
\end{align*}
for all smooth functions $\phi$, where ``h.o.t.'' stands for
``higher order terms''. 
\end{lemma}
\begin{remark}
Due to the convexity of the equation,
in this example the upper and lower bounds are different.
\end{remark}
\begin{remark}
We have only stated the principal error terms, the terms deciding the
rate. The other terms are put in the ``h.o.t.'' category. Since the
principal error terms need not be the lowest order terms (see the first
inequality in Lemma \ref{ssLem1}), maybe a better name than
``h.o.t.'' would be the ``less important terms''.
\end{remark}
A direct consequence of Proposition \ref{semigr_main} is the following
result:  
\begin{corollary}
\label{ssCor1}
Let $u_h$ denote the solution of \eqref{Splitt} where $S_h$ is defined
in \eqref{SSplitt} and $u_{h,0}=u_0$, and let $u$ be the solution of
\eqref{spliteq00} with initial value $u_0$. Under the assumptions of
Lemma \ref{ssLem1} we have  
$$-C\dt^{\frac1{13}}\leq u-u_h\leq
C\dt^{\frac29}\quad \text{in}\quad \dt\{0,1,2,\dots,n_T\}\times\R^N.$$
\end{corollary}

Next, we give the result when $F_1$ and $F_2$ are assumed to be only Lipschitz
continuous (which is the natural regularity assumption here). In this
case the consistency relation is:
\begin{lemma} 
\label{ssLem2}
Under the above assumptions, if $F_1$ and $F_2$ are only Lipschitz
continuous, we have
\begin{align*}
&\Big|\frac1{\dt}[S_h(t)-1]\phi(t_{n-1},x)+F_1(D^2\phi(t_n,x)) +
  F_2(D^2\phi(t_n,x))\Big|\\ 
&\leq  C\dt |D^2\phi_t|_0 + C\dt^{\frac12}|D^3\phi|_0  + \text{h.o.t.}
\end{align*}
for all smooth functions $\phi$.
\end{lemma}

Again as a direct consequence of Proposition \ref{semigr_main} we have
the following error bound: 
\begin{corollary}
\label{ssCor2}
Under the assumptions of Corollary \ref{ssCor1} but where $F_1$ and
$F_2$ are only assumed to be Lipschitz continuous, we have 
\begin{align*}
-C\dt^{\frac1{14}}\leq u-u_h\leq
  C\dt^{\frac1{6}}\quad \text{in}\quad \dt\{0,1,2,\dots,n_T\}\times\R^N.
\end{align*}
\end{corollary}

\begin{remark}
We see a slight reduction of the rates in the Lipschitz case but not
as important as one might have guessed.  For first order equations
these methods lead to the same rates in the smooth and Lipschitz
cases.
\end{remark}

\begin{remark}
If we change operators $S_1,S_2$ so that $S_1(t)\phi$ and
$S_2(t)\phi$ denote the viscosity solutions of
\begin{align*}
&u(x)=\phi(x)-t F_1(D^2u(x)) \quad \text{in}\quad\R^N,\\
&u(x)=\phi(x)-t F_2(D^2u(x)) \quad \text{in}\quad\R^N,
\end{align*}
respectively, then the statements of Corollary \ref{ssCor1} and
\ref{ssCor2} still hold. 
\end{remark}

In the proofs of Lemmas \ref{ssLem1} and \ref{ssLem2} we will use the
 following lemma:
\begin{lemma}
\label{ssLem3}
Let $\bar S$ be the semigroup associated to the equation 
$$u_t+\bar F(D^2u)=0,$$
where  $\bar F$ is Lipschitz, convex, and non-increasing. Define $\bar
F_\delta$ by 
$$\bar F_\delta=\bar F*\bar\rho_\delta,$$
where $\bar\rho_\delta(X)=\delta^{-N^2}\bar \rho(X/\delta)$ and $\bar
\rho$ is a smooth function on $\BS(N)$ with mass one and support in $B(0,1)$.
Then for any smooth function $\phi$,
$$\bar S(t)\phi-\phi+t\bar F_\delta(D^2\phi)\leq t\delta |D\bar F|_0+\frac12
t^2|D\bar F|_0|D\bar F_\delta|_0|D^4\phi|_0,$$
and 
$$\bar S(t)\phi-\phi+t\bar F_\delta(D^2\phi)\geq -\frac12
t^2|D\bar F|_0(|D^2\bar F_\delta|_0|D^3\phi|_0^2+|D\bar
F_\delta|_0|D^4\phi|_0).$$  
\end{lemma}
The proof of this result will be given after the proofs of Lemmas
\ref{ssLem1} and \ref{ssLem2}.

\begin{proof}[Proofs of Lemmas \ref{ssLem1} and \ref{ssLem2}]
In order to treat the two results at the same time, we mollify $F_1$
and $F_2$ and consider $F_{1,\delta}$ and $F_{2,\delta}$ (see Lemma
\ref{ssLem3} for the definitions). 
By Lemma \ref{ssLem3} we have the following (small time) expansions:
\begin{align}
\label{ST1}
&S_j(t)\phi-\phi+t F_{j,\delta}(D^2\phi)\leq t\delta |D F_j|_0+\frac12
t^2|D F_j|_0|D F_{j,\delta}|_0|D^4\phi|_0,\\
\label{ST2}
&S_j(t)\phi-\phi+t F_{j,\delta}(D^2\phi)\\
&\nonumber
\geq -\frac12
t^2|D F_j|_0(|D^2 F_{j,\delta}|_0|D^3\phi|_0^2+|DF_{j,\delta}|_0|D^4\phi|_0),
\end{align}
for smooth functions $\phi$ and $j=1,2$.

Now we want to find an (small time) expansion for $S_h$. We write
\begin{align*}
&S_h(t)\phi-\phi+t(F_1+F_2)(D^2\phi)\\
&= [S_1(t)S_2(t)\phi-S_1(t)(\phi-tF_{2,\delta}(D^2\phi))]\\
&\quad +[S_1(t)(\phi-tF_{2,\delta}(D^2\phi))-\phi
  +tF_{1,\delta}(D^2\phi)+tF_{2,\delta}(D^2\phi)]\\
&\quad + t[(F_1+F_2)(D^2\phi)-(F_{1,\delta}+F_{2,\delta})(D^2\phi)].
\end{align*}
In view of the Lipschitz regularity and convexity of $F_1$ and $F_2$,
the last term on right hand side is between $-Ct\delta$ (Lipschitz
regularity) and $0$ (convexity),    
while the first 2 terms can be estimated using non-expansiveness and
small time expansions for $S_1$ and $S_2$. The principal error term comes from
the small time expansion for the term $S_1(t)(\phi-tF_{2,\delta}(D^2\phi))$.
In view of \eqref{ST1} and \eqref{ST2},
$$t\delta |DF_1|_0+\frac12
t^2|D F_1|_0|D F_{1,\delta}|_0|D^4\{\phi-tF_2^\delta(D^2\phi)\}|_0$$
is an upper bound on the principal error term, while
$$\frac12t^2|D
F_1|_0\Big[|D^2F_{1,\delta}|_0|D^3\{\phi-tF_{2,\delta}(D^2\phi)\}|_0^2
+|DF_{1,\delta}|_0|D^4\{\phi-tF_{2,\delta}(D^2\phi)\}|_0\Big]$$
is a lower bound. Expanding out these expressions keeping only the ``worst
terms'' and bearing in mind the Lipschitz regularity of $F_1$ and $F_2$,
lead to the following upper 
and lower bounds respectively,
$$C(t\delta+t^2|D^4\phi|_0 + t^3|D^4F_{2,\delta}|_0|D^3\phi|_0^4)
\quad\text{and}$$
$$C|D^2F_{1,\delta}|_0\Big(t^2|D^3\phi|^2 +
t^3|D^3F_{2,\delta}|_0|D^3\phi|_0^4  
+t^4 |D^3F_{2,\delta}|^2|D^3\phi|_0^6\Big).$$  

To conclude the proofs of the upper bounds in Lemmas \ref{ssLem1} and
\ref{ssLem2}, we note that 
\begin{align*}
&\frac1{\dt}[S_h(t)-1]\phi(t_{n-1},x)+F_1(D^2\phi(t_n,x)) +
  F_2(D^2\phi(t_n,x))\\ 
&\leq \dt (|DF_1|_0+|DF_2|_0)|D^2\phi_t|_0 \\
&\quad +
  \Big[\frac1{\dt}[S_h(t)-1]\phi+F_1(D^2\phi)+
  F_2(D^2\phi)\Big]_{(t_{n-1},x)}.
\end{align*}
In view of the above estimates the right hand side can be upper bounded by
\begin{align}
\label{ssExpr}
C\Big[ \dt|D^2\phi_t|_0 +\delta +
  \dt|D^4\phi|_0+\dt^2|D^4F_{2,\delta}|_0|D^3\phi|_0^4\Big]. 
\end{align}
This proves the upper bound in Lemma \ref{ssLem1} after sending
$\delta\ra 0$ while keeping in mind that in this case,
$$|D^nF^\delta_2|_0\leq
|D^nF_2|_0<\infty\quad\text{and}\quad|D^mF^\delta_1|_0\leq 
|D^mF_1|_0<\infty$$
for $n=1,\dots,4$ and $m=1,2$. To get the upper bound in Lemma
\ref{ssLem2}, we only need to note that in this case $|D^nF_j|_0\leq C
\delta^{1-\delta}$, $n\in\N$, $j=1,2$, and then minimize \eqref{ssExpr}
w.r.t. $\delta$. 

The upper bounds follow in a similar way. We conclude the proof simply
by giving the expression corresponding to \eqref{ssExpr},
\begin{align*}
&C\Big[ \dt|D^2\phi_t|_0 +\delta \\
&+  |D^2F_{1,\delta}|_0\Big(\dt|D^3\phi|^2 +
  \dt^2|D^3F_{2,\delta}|_0|D^3\phi|_0^4 +\dt^3
  |D^3F_{2,\delta}|^2|D^3\phi|_0^6\Big)\Big].
\end{align*}


\end{proof}

\begin{proof}[Proof of Lemma \ref{ssLem3}]
Let 
$$w=\phi-t\bar F_\delta(D^2\phi),$$
and observe that
\begin{align}
\label{PFeq1}
w_t+\bar F(D^2w)=-\bar F_\delta(D^2\phi)+\bar F(D^2\phi)-\bar
F(D^2\phi)+\bar F(D^2w). 
\end{align}
Since $\bar F$ is convex, it is easy to see that $\bar F_\delta(X)\geq
\bar F(X)$, and hence 
\begin{align}
\label{PFeq2}
-|D\bar F|_0\delta \leq -\bar F_\delta(D^2\phi)+\bar F(D^2\phi)\leq 0.
\end{align}
The second difference, $-\bar F(D^2\phi)+\bar F(D^2w)$, can be written
\begin{align}
&\int^1_0\frac d{ds}\{\bar F(sD^2w+(1-s)D^2\phi)\}ds\nonumber\\
&= \sum_{ij} \del_i\del_j(w-\phi)\int_0^1
(\del_{X_{ij}} \bar F)(sD^2w+(1-s)D^2\phi)ds\nonumber\\
&=-t \sum_{ij} \del_i\del_j\{\bar F_\delta(D^2\phi)\}\int_0^1
(\del_{X_{ij}} \bar F)(sD^2w+(1-s)D^2\phi)ds.\label{PFeq3}
\end{align}
We expand $\del_i\del_j\{\bar F_\delta(D^2\phi)\}$ and get
\begin{align*}
&\sum_{klmn}(\del_{X_{kl}}\del_{X_{mn}}\bar F_\delta)(D^2\phi)(\del_i\del_k\del_l\phi)(\del_j\del_m\del_n\phi)\\
&+\sum_{kl}(\del_{X_{kl}}\bar F_\delta)(D^2\phi)(\del_i\del_j\del_k\del_l\phi).
\end{align*}
We call the first term $M[\phi]_{ij}$. 

Since $\del_{X_{kl}}\del_{X_{mn}}\bar
F_\delta=\del_{X_{mn}}\del_{X_{kl}}\bar F_\delta$, 
it follows that $M$ is symmetric,
$$M[\phi]_{ij}=M[\phi]_{ji}.$$
Moreover, since $\bar F$ is convex, $M$ is
positive semidefinite: For every $\xi\in\R^N$
\begin{align*}
&\sum_iM[\phi]_{ij}\xi_i\xi_j\\
&=\sum_{klmn}(\del_{X_{kl}}\del_{X_{mn}}\bar F_\delta)(D^2\phi)(\del_k\del_l(\sum_i\xi_i\del_i\phi))(\del_m\del_n(\sum_j
\xi_j\del_j\phi))\\
&=\sum_{klmn}(\del_{X_{kl}}\del_{X_{mn}}\bar F_\delta)(D^2\phi)Y_{kl}Y_{mn}
\geq 0,
\end{align*}
where $Y_{ij}=\del_i\del_j(\sum_k\xi_k\del_k\phi)$ and where the
inequality follows by convexity of $\bar F$.

By the spectral theorem there exists $e^k\in\R^N$ and $\lambda_k\in\R$
for $k=1,\dots,N$ (depending on $\phi$) such that 
$$M[\phi]=\sum_k\lambda_k e^k\otimes e^k.$$
Furthermore, since $M$ is positive semidefinite, $\lambda_i\geq0$ for
$i=1,\dots,N$. Therefore we have 
\begin{align*}
&\sum_{ij}M[\phi]_{ij}\int_0^1
(\del_{X_{ij}} \bar F)(sD^2w+(1-s)D^2\phi)ds\\
&=\sum_k \lambda_k \int_0^1\sum_{ij}e^k_ie^k_j
(\del_{X_{ij}} \bar F)(sD^2w+(1-s)D^2\phi)ds \leq 0,
\end{align*}
where the inequality follows from the fact that $\bar F$ is
nonincreasing. We conclude that 
$$\bar F(D^2w)-\bar F(D^2\phi)\geq -t|D\bar F|_0|D^2\bar F_\delta|_0|D^4\phi|_0
,$$
and hence by \eqref{PFeq1} -- \eqref{PFeq3} we get
$$w_t+\bar F(D^2w)\geq
-|D\bar F|_0\delta-t|D\bar F|_0|D^2\bar F_\delta|_0|D^4\phi|_0.$$
The first part of the Lemma now follows from the comparison principle.

The second part of the Lemma follows from \eqref{PFeq1} --
\eqref{PFeq3} and the comparison principle after noting that this
time, due the its sign, the $D^2\bar F_\delta$ term will be part of the
error expression. 
\end{proof}

\subsection{Piecewise constant controls}
Here we study approximations by piecewise constant controls. Such
approximations have been studied e.g. in \cite{LM,Kr:Const} (see also
the references therein). We consider the following
simplified version of equation \eqref{E},
\begin{align}
\label{KS1}
u_t+\max_{i}\{-L^i u -f^i(x)\}=0 \qquad \text{in}\quad Q_T,
\end{align}
where
$$L^i\phi= \tr[\sigma^i(x)\sigma^{i\,T}(x)D^2\phi]+b^i(x)D\phi+c^i(x)\phi,$$
and $\sigma,b,c$ and $f$ satisfy assumption (A1) when $\alpha$ is replaced
by $i$. Note that the coefficients are independent of time. We
approximate \eqref{KS1} in the following way,
\begin{align}
u^{n+1}(x)=\min_i S_i(\dt)u^n(x) \qquad \text{in} \quad
\{0,1,\dots,n_T\}\times\R^N, 
\end{align}
where $S_i(t)\phi(x)$ denotes the solution at $(t,x)$ of the linear equation
\begin{align}
\label{KSLin}
u_t-L^i u - f^i(x) = 0
\end{align}
with initial data $\phi$ at time $t=0$. As usual, $u^n$ is expected to
be an approximation of $u(t_n,x)$, $t_n:=n\dt$, and we are looking
for a bound on the approximation error.

Under assumption (A1) the comparison principle holds for the linear
equations \eqref{KSLin}, hence $S_i$ and $\min_iS_i$ are
monotone. Furthermore, we have the following consistency relation:

\begin{lemma}
\label{KSLem}
If (A1) holds, then for any smooth function $\phi$ we have
\begin{align*}
&\Big|\frac1{\dt}[\min_iS_i(\dt)-1]\phi(t_{n-1},x) +
  \dt\max_i\{-L^i\phi(t_n,x) - f^i(x)\}|\\ 
&\leq C\dt|D^2\phi_t|_0+C\dt\left(\sum_{n=0}^4|D^n\phi|_0+1\right) +
  C\dt^{1/2}\left(\sum_{n=0}^2|D^n\phi|_0+1\right).
\end{align*}
\end{lemma}

We have the following error bound: 

\begin{proposition}
\label{KSCor}
Assume (A1). Let $u_h$ denote the solution of \eqref{Splitt}
corresponding to 
$$S_h(t_{n-1},t_n)=\min_iS_i(\dt)$$ 
and $u_{h,0}=u_0$, and let $u$ be the solution of \eqref{KS1} with
initial value $u_0$. Then
$$-C\dt^{\frac1{10}}\leq u-u_h\leq
0 \quad \text{in}\quad \dt\{0,1,2,\dots,n_T\}\times\R^N.$$
\end{proposition}

\begin{proof}
We first observe that $u_h
\geq u$ in $Q_T$. This can be easily seen from the control
  interpretation of $u_h$ (which we have not provided!) or from the
  comparison principle since $u_h$ is a supersolution of \eqref{KS1}
  (solutions of \eqref{KSLin} are supersolutions of \eqref{KS1} and so
  is the min of such solutions). The other bound follows from Lemma
  \ref{KSLem} and Proposition \ref{semigr_main}.
\end{proof}
\begin{remark}
Assuming more regularity on the coefficients does not lead to any
improvement of the bound. The principal contribution to the error
comes from the $|D^4\phi|_0$-term, and this term does not depend on
the regularity of the coefficients (only on the $L^\infty$ norm of $\sigma$). 
\end{remark}
\begin{remark}
\label{Rem:Const}
In \cite{Kr:Const} Krylov obtains a better rate, namely $1/6$. His
approach is different for ours, he works on the
dynamic programming principle directly using control techniques.
\end{remark}

\begin{proof}[Proof of Lemma \ref{KSLem}] 

Let $\sigma^i_\delta=\sigma^i*\rho_\delta$,
and define similarly $b^i_\delta$, $c^i_\delta$, and $f^i_\delta$, and
let $L^i_\delta$ be the operator $L^i$ corresponding to
$\sigma^i_\delta,b^i_\delta,c_\delta^i$. Observe that 
\begin{align*}
&|\min_iS_i(t)\phi -\phi + t\max_i\{-L_\delta^i\phi - f^i_\delta(x)\}|\\
&=|\min_i(S_i(t)\phi -\phi)-t\min_i\{L_\delta^i\phi + f^i_\delta(x)\}|\\
&\leq  \max_i|S_i(t)\phi -\phi + t(-L_\delta^i\phi - f^i_\delta(x)) |.
\end{align*}
Next, define
$$w^\pm=\phi - t(-L_\delta^i\phi -
f^i_\delta(x))\pm\frac12t^2|L^iL_\delta^i\phi-L^if^i_\delta|_0 \pm
t|(L_\delta^i-L^i)\phi - (f^i_\delta-f^i)|_0,$$
and observe that $w^+$ is a supersolution of \eqref{KSLin} while $w^-$
is a subsolution. By the comparison principle and
properties of mollifiers we get
\begin{align*}
&|S_i(t)\phi -\phi + t(-L_\delta^i\phi - f^i_\delta(x))|\\
&\leq \frac12 t^2
|L^iL_\delta^i\phi-L^if^i_\delta|_0+t\delta
C\left(\sum_{n=0}^2|D^n\phi|_0+1\right).
\end{align*} 
Furthermore, by properties of mollifiers and the Lipschitz regularity
of the coefficients we see that
$$|L^iL_\delta^i\phi+L_if_i|_0\leq
C\left(\sum_{n=0}^4|D^n\phi|_0 +\delta^{-1}\sum_{n=0}^2|D^n\phi|_0
+\delta^{-1}+1\right).$$  
By combining the above estimates we get
\begin{align*}
&|\min_iS_i(t)\phi -\phi + t\max\{-L^i\phi - f^i(x)\}|\\
&\leq Ct^2\left(\sum_{n=0}^4|D^n\phi|_0 +
  \delta^{-1}\sum_{n=0}^2|D^n\phi|_0 
+\delta^{-1}+1\right) + t\delta
C\left(\sum_{n=0}^2|D^n\phi|_0+1\right), 
\end{align*}
and the result follows by similar arguments as was given in the proofs
of Lemmas \ref{ssLem1} and \ref{ssLem2} after optimizing
w.r.t. $\delta$. 
\end{proof}

\section{Remarks on the H{\"o}lder continuous case}
\label{Sec:Rem}

In this section we give an extension of the main
result Theorem \ref{mainres} to the case when solutions of \eqref{E}
do no longer belong to the space $\CC^{0,1}$ but rather belong to the
bigger space $\CC^\beta$ for some $\beta\in(0,1)$.

In the time-dependent case $\CC^\beta$ regularity of the solution is
observed typically when assumption (A1) is relaxed in the following way:
\medskip

\noindent {\bf (A1')} For any $\alp \in \A$,
$a^{\alp}=\frac12\sigma^{\alp}{\sigma^{\alp}}^T$ for some $N\times P$ 
matrix $\sigma^{\alp}$. Moreover, there is a constant $K$ independent
of $\alp$ such that
$$|u_0|_\beta+|\sigma^{\alp}|_1+|b^{\alp}|_1 +
|c^{\alp}|_\beta+|f^{\alp}|_\beta \leq 
K.$$

In other words $u_0, c^{\alp}, f^{\alp}$ now belongs to
$\CC^\beta$. 

\begin{lemma}
If (A1') holds, then there exists a unique solution $u\in
\CC^{0,\beta}(\ol Q_T)$ of \eqref{E} and \eqref{IV}.
\end{lemma}

This standard result is proved e.g. in \cite{JK:ContDep}. We claim
that under (A1'), we have the same regularity (the same $\beta$) for
all equations considered in this paper. We skip the proof of this claim. In
the rest of this section, the solutions of the different equations
belong to $\CC^{0,\beta}(\ol Q_T)$ with the same fixed $\beta\in(0,1]$. 

Lower than $\CC^{0,1}$ regularity of solutions implies lower
convergence rates than obtained in Sections \ref{Sec:Sw} --
\ref{Sec:Appl}. We will now state the H{\"o}lder versions of some these
results without proofs. The proofs are not much different from the
proofs given above, and moreover, the H{\"o}lder case was extensively
studied in \cite{BJ:Rate}. We start by the convergence rate for the
switching system approximation of Section~\ref{Sec:Sw}.
\begin{proposition}
\label{sw-rate2}
Assume (A1'). If $\bar{u}$ and $v$ are the solutions of \eqref{HJBi} and
\eqref{Sw1} in $\CC^{0,\beta}(\ol Q_T)$, then for $k$ small enough, 
$$0\leq v_i -\bar{u}\leq Ck^{\frac{\beta}{2+\beta}}\quad
\text{in}\quad\ol Q_T,\quad i\in\II,$$ 
where $C$ only depends on $T$ and $K$ from (A1').
\end{proposition}

In order to state a $\CC^\beta$ version of Theorem \ref{mainres} we
need to modify assumption (S3). The requirement on $\phi_\eps$ should
be changed to
$$|\partial_t^{\beta_0}D^{\beta'} \phi_\e (x,t) | \leq \tilde K
\e^{\beta-2\beta_0-|\beta'|}  \quad \hbox{in  }\overline Q_T ,\quad\text{ for
any $\beta_0\in\N$, $\beta'\in\N^{N}$}.$$
We will denote the modified assumption by (S3'). Now we state the
$\CC^\beta$ version of our main result, Theorem \ref{mainres}.

\begin{theorem}
Assume (A1'), (S1), (S2) and that \eqref{S} has a unique
solution $u_h\in C_b(\G_h)$. Let $u$ denote the solution of
\eqref{E}-\eqref{IV}, and let $h$ be sufficiently small. 

(a) {\bf (Upper bound)} If (S3')(i) holds, then there exists a constant
$C$ depending only $\mu$, $K$ in (S1), (A1') such that
$$ u-u_h \leq e^{\mu t}|(u_0-u_{0,h})^+|_0 + C\min_{\e>0} \left(\e^\beta
+ E_1 (\tilde K,h,\e)\right) \quad\text{in}\quad\G_h,
$$
where $\tilde K = |u|_1$.

(b) {\bf (Lower bound)} If (S3')(ii) and (A2) holds, then there exists
a constant $C$ depending only $\mu$, $K$ in (S1), (A1') such that
$$
u-u_h \geq -e^{\mu t}|(u_0-u_{0,h})^-|_0 - C\min_{\e>0}
\left(\e^{\frac{\beta^2}{2+\beta}} + E_2 
(\tilde K,h,\e)\right) \quad\text{in}\quad\G_h, 
$$
where $\tilde K = |u|_1$.
\end{theorem}

\begin{remark}
For the FDMs described in Section \ref{Sec:FDM} we get an upper rate of
$\frac{\beta}2$ and a lower rate of $\frac{2\beta^2}{4(2+\beta)-2\beta}$ in
the $\CC^\beta$ case. Compare with Theorem \ref{rate1}.
\end{remark}

\appendix
\section{Well-posedness, regularity, and continuous dependence for
switching systems}
In this section we give well-posedness, regularity, and
continuous dependence results for solutions of a very general
switching system that has as special cases the scalar HJB equations
\eqref{E}, and the switching systems \eqref{Sw1}, \eqref{Sw2},
\eqref{Sw3}. 

We consider the following system:
\begin{gather}
\label{Sw}
F_i(x,u,\del_tu_i,Du_i,D^2u_i)=0 \quad \text{in}\quad Q_T, \quad
i\in\II:=\{1,\dots,M\}, 
\end{gather}
with 
\begin{align*}
F_i(t,x,r,p_t,p_x,X)&
=\max\Big\{p_t +
\sup_{\alp\in\A}\inf_{\beta\in\B}\LL^{\alp,\beta}_i(x,r_i,p_x,X);\   
r_i-\M_ir\Big\},\\
\LL^{\alp,\beta}_i(t,x,s,q,X)&=-\tr[a^{\alp,\beta}_i(t,x) X] -
b^{\alp,\beta}_i(t,x) q - c^{\alp,\beta}_i(t,x) s -f^{\alp,\beta}_i(t,x),
\end{align*}
where $\M$ is defined below \eqref{Sw1}, $\A,\B$ are compact metric
spaces, $r$ is a vector $r=(r_1,\dots,r_M)$, and $k>0$ is a constant (the 
switching cost). See \cite{EF:Sw,CDE:Sw,Ya:Sw2,IK:Sw,IK:Sys} for
more information about such systems. 

We make the following assumption:\\

\noindent {\bf (A)} For any $\alp,\beta,i$, 
$a^{\alp,\beta}_i=\frac12\sigma^{\alp,\beta}_i{\sigma^{\alp,\beta}_i}^T$
for some $N\times P$ matrix $\sigma^{\alp,\beta}_i$. Furthermore,
there is a constant $C$ 
independent of $i,\alp,\beta,t$, such that 
$$|\sigma^{\alp,\beta}_i(t,\cdot)|_1+|b^{\alp,\beta}_i(t,\cdot)|_1 
+|c^{\alp,\beta}_i(t,\cdot)|_1+|f^{\alp,\beta}_i(t,\cdot)|_1\leq \bar
C.$$\\[-0.5cm]

We start by comparison, existence, uniqueness, and $L^\infty$ bounds on the
solution and its gradient. Before stating the results, we define
$USC(\bar Q_T;\R^M)$ and $LSC(\bar Q_T;\R^M)$ to be the spaces of
upper and lower semi-continuous functions from $\bar Q_T$ into $\R^M$
respectively. 

\begin{theorem}
\label{WP}
Assume (A) holds. 

(i) If $u\in USC(\bar Q_T;\R^M)$ is a subsolution of \eqref{Sw} bounded above
and $v\in LSC(\bar Q_T;\R^M)$ supersolution of \eqref{Sw} bounded below, then
$u\leq v$ in $\bar Q_T$.

(ii) There exists a unique bounded continuous solution $u$ of \eqref{Sw}. 

(iii) The solution $u$ of \eqref{Sw} belongs to $\CC^{0,1}(\bar Q_T)$,
and satisfies for all $t,s\in[0,T]$
\begin{align*}
e^{-\lambda t}\max_i |u_i(t,\cdot)|_0\leq \ & \max_i
  |u_{0,i}|_0+t\sup_{i,\alp,\beta} |f^{\alp,\beta}_i|_0,\\
\intertext{where $\lambda:=\sup_{i,\alp,\beta}|c^{\alp,\beta+}_i|_0$,}
e^{\lambda_0 t}\max_i\,[u_i(t,\cdot)]_1\leq \ & 
  \max_i[u_{0,i}]_1+t
  \sup_{i,\alp,\beta,s} 
\Big\{|u^i|_0[c^{\alp,\beta}_i(s,\cdot)]_1+[f^{\alp,\beta}_i(s,\cdot)]_1\Big\},
\end{align*}
where
$\lambda_0:=\sup_{i,\alp,\beta,s}\{|c^{\alp,\beta+}_i(s,\cdot)|_0
+[\sigma^{\alp,\beta}_i(s,\cdot)]_1^2+[b^{\alp,\beta}_i(s,\cdot)]_1\}$, 
and 
\begin{align*}
\max_i|u_i(t,x)-u_i(s,x)|\leq \ &C|t-s|^{1/2},
\end{align*}
where $C\leq 8 M \bar C+\sqrt T  \bar C (2M+1)$ and
$M:=\sup_{i,t}|u_i(t,\cdot)|_1$.
\end{theorem}

Before giving the proof we state a key technical lemma.
\begin{lemma}
\label{SYS2SC}
Let $u\in USC(\bar Q_T;\R^M)$ be a bounded above subsolution of
\eqref{Sw} and
$\bar{u}\in LSC(\bar Q_T;\R^M)$ be a bounded below supersolution
of an other equation \eqref{Sw} where the functions
$\LL^{\alp,\beta}_i$ are replaced by functions
$\bar{\LL}^{\alp,\beta}_i$ satisfying the same
assumptions. Let $\phi:[0,T]\times\R^{2N}\ra \R$ be a smooth function
bounded from below. We denote by 
$$
\psi_i(t,x,y)=u_i(t,x)-\bar{u}_i(t,y)-\phi(t,x,y)\; ,$$ 
and $M=\sup_{i,t,x,y}\,\psi_i(t,x,y)$. 
If there exists a maximum point for $M$, i.e. a point $(i',t_0,x_0,y_0)$
such that $\psi_{i'}(t_0,x_0,y_0) = M$, then there exists $i_0 \in \II$
such that $(i_0,t_0,x_0,y_0)$ is also a maximum point for $M$, and, in
addition $\bar{u}_{i_0}(t_0,y_0)<\M_{i_0}\bar{u}(t_0,y_0)$.
\end{lemma}

Loosely speaking this lemma means that whenever we do doubling of
variables for systems of the type \eqref{Sw}, we can ignore the
$u_i-\M_iu$ parts of the equations.  So we are more or less back in the
scalar case with equations
$\del_tu_{i_0}+\sup_\alp\inf_\beta\LL_{i_0}^{\alp,\beta}[u_{i_0}]\leq 0$ and
$\del_t\bar{u}_{i_0}+
\sup_\alp\inf_\beta\bar{\LL}^{\alp,\beta}_{i_0}[\bar{u}_{i_0}]\geq0$.
We skip the proof since it is similar to the proof given in \cite{BJ:Rate2}
for the stationary case. 





\begin{proof}[Proof of Theorem \ref{WP}]
Comparison, uniqueness, and existence is proved in \cite{IK:Sw} for
the stationary Dirichlet problem for \eqref{E} on a bounded domain under
similar assumptions on the data. To extend the comparison result to
a time dependent problem in an unbounded domain, we only need to
modify the test function used in \cite{IK:Sw} in the standard
way. (See also the arguments given below). Comparison implies
uniqueness, and existence follows from Perron's method. This last
argument is similar to the argument given in \cite{IK:Sw}, but easier
since we have no boundary conditions other than the initial condition.

Let 
$$w(t):= e^{\lambda t}\Big\{\max_i
  |u_{0,i}|_0+t\sup_{i,\alp,\beta} |f^{\alp,\beta}_i|_0\Big\},\quad 
$$
then the bound on $|u|_0$ follows from the comparison principle
after checking that $w$ ($-w$) is a supersolution (subsolution) of
\eqref{Sw}. 

To get the bound on the gradient of $u$, consider
$$m:=\sup_{i,t,x,y\in\R^N}\left\{u_i(t,x)-u_i(t,y)-\bar
w(t)|x-y|\right\},$$
where
\begin{align*} 
\bar w(t):=\ & e^{\lambda_0 t}\Big\{\max_i[u_{0,i}]_1+t
  \sup_{i,\alp,\beta,s} 
\Big\{|u_i|_0[c^{\alp,\beta}_i(s,\cdot)]_1 +
  [f^{\alp,\beta}_i(s,\cdot)]_1\Big\}\Big\}. 
\end{align*} 
We are done if we can prove that $m\leq 0$. Assume this is not the
case, $m>0$, and for simplicity that this maximum is attained in
$\bt,\bx,\by$. Then there exists a $k>0$ such that 
$$u_i(\bt,\bx)-u_i(\bt,\by)-\bar w(\bt)|\bx-\by| - \bt e^{\lambda_0
  \bt}k>0, \quad i\in\II.$$
Let $\psi_i(t,x,y):=u_i(t,x)-u_i(t,y)-\bar w(t)|x-y| - t e^{\lambda_0
  t}k$, then $\psi$ also has maximum $M>0$ at some point
  $(\tilde i, \tilde t, \tilde x,\tilde y)$. Since $M>0$,
  $\tilde x\neq\tilde y$ and $\tilde t>0$. Therefore $\bar
  w(t)|x-y|+te^{\lambda_0 t}k$ is a smooth function at $(\tilde t,
  \tilde x,\tilde y)$ and a standard argument using
the viscosity sub- and supersolution inequalities for \eqref{Sw} at
  $(\tilde t, \tilde x,\tilde y)$ and Lemma \ref{SYS2SC} leads to
  $k\leq 0$. See the proof of  Theorem \ref{CD} for a similar argument.
This is a contradiction and hence $m\leq 0$.

In the general case when the maximum $m$ need not be attained at some
finite point, we must modify the test function in the standard way. We
skip the details.

To get the time regularity result, assume that $s<t$ and let $u^\eps$
be the solution of \eqref{Sw} in $t\in(s,T]$ starting from
  $u(s,\cdot)*\rho_\eps(x)=:u^\eps_0(x)$. By the comparison principle
$$|u-u^\eps|\leq \sup_{r\in[s,T]}[u(r,\cdot)]_1\eps \quad
  \text{in}\quad [s,T]\times\R^N,$$
and easy computations show that
$$w^\pm_\eps(t,x)=e^{\lambda t}\Big\{u^\eps_0(x)\pm (t-s) C_\eps \Big\}$$ 
are subsolution ($w^-$) and supersolution ($w^+$) of \eqref{Sw} if
$$C_\eps=\bar C^2|D^2u^\eps_0|_0+\bar C(|Du^\eps_0|_0+|u^\eps_0|_0+1)$$
and $\bar C$ is given by (A). 
Another application of the comparison principle then yields
$$w^-_\eps\leq u^\eps \leq w^+_\eps\quad
  \text{in}\quad [s,T]\times\R^N.$$
The result now follows from 
\begin{align*}
&|u(t,x)-u(s,x)| \\
& \leq
  |u(t,x)-u^\eps(t,x)|+|u^\eps(t,x)-u^\eps_0(x)|+|u^\eps_0(x)-u(s,x)|\\
&\leq ([u(t,\cdot)]_1+[u(s,\cdot)]_1)\eps + |t-s|C_\eps,
\end{align*}
and a minimization in $\eps$ after noting that $C_\eps\leq C(\eps^{-1}+1)$.
\end{proof}

We proceed to obtain continuous dependence on the coefficients.
\begin{theorem}
\label{CD}
Let $u$ and $\bu$ be solutions of \eqref{Sw} with coefficients
$\sigma,b,c,f$ and $\bs,\bb,\bc,\bar{f}$ respectively. If both sets of
coefficients satisfy (A1), and
$|u|_0+|\bu|_0+[u(t,\cdot)]_1+[\bu(t,\cdot)]_1\leq M<\infty$ for $t\in
[0,T]$, then 
\begin{align*}
&e^{-\lambda t}\max_i|u_i(t,\cdot)-\bu_i(t,\cdot)|_0\leq
\max_i|u_{i}(0,\cdot)-\bu_i(0,\cdot)|_0\\ 
&+t^{1/2} K \sup_{i,\alp,\beta}|\sigma-\bs|_0
+t \sup_{i,\alp,\beta}\Big\{2M|b-\bb|_0+M|c-\bc|_0+|f-\bff|_0\Big\},
\end{align*} 
where $\lambda:=\sup_{i,\alp,\beta}|c^-|_0$ and
\begin{align*}
K^2\leq &\ 8M^2 +
8MT\sup_{i,\alp,\beta}\Big\{2M[\sigma]^2_1\wedge[\bs]^2_1\\ 
&+2M[b]_1\wedge[\bb]_1+M[c]_1\vee[\bc]_1+[f]_1\wedge[\bff]_1\Big\}. 
\end{align*} 
\end{theorem}

\begin{proof}
We only indicate the proof in the case $\lambda=0$. Define
\begin{align*}
&\psi^i(t,x,y):=u_i(t,x)-\bu_i(t,y)-\frac{1}{\delta}|x-y|^2-\eps(|x|^2+|y|^2),\\
&m:=\sup_{i,t,x,y}\psi^i(t,x,y)-\sup_{i,x,y}(\psi^i(0,x,y))^+,\\
&\bar m:=\sup_{i,t,x,y}\left\{\psi^i(t,x,y)-
  \frac{\sigma m t}T\right\},
\end{align*}
where $\sigma\in(0,1)$. We assume $m>0$ since otherwise we are
done. We will now derive an upper bound on $m$. To do this we consider
$\bar m$. By the assumptions this supremum is attained at some point
$(i_0,t_0,x_0,y_0)$. Since $m>0$ it follows that $\bar m>0$ and $t_0>0$,
and by Lemma \ref{SYS2SC}, the index $i_0$ may be chosen  
so that $\bar{u}_{i_0}(t_0,y_0)<\M_{i_0}\bar{u}(t_0,y_0)$. With this in mind,
the maximum principle for semi continuous functions
\cite{CI:MaxPr,CIL:UG} and the definition of viscosity solutions imply
the following inequality:
$$p_t-\bar
p_t+\sup_\alp\inf_\beta\LL^{\alp,\beta}_{i_0}(t_0,x_0,u_{i_0},p_x,X)-\sup_\alp\inf_\beta\bar{\LL}^{\alp,\beta}_{i_0}(t_0,y_0,\bu_{i_0},p_y,Y)\leq0,$$ 
where $(p_t,p_x,X)\in \ol{\mathcal P}^{2,+}u_{i_{0}}(x_0)$ and $(\bar
p_t, p_y,Y)\in \ol{\mathcal P}^{2,-}\bu_{i_{0}} (y_0)$ (see
\cite{CI:MaxPr,CIL:UG} for the notation). Furthermore
$p_t-\bar p_t= \frac{\sigma m}{T}$,
$p_x=\frac2\delta(x_0-y_0)+2\eps x_0$, 
$p_y=\frac2\delta(x_0-y_0)-2\eps y_0$, and  
$$\begin{pmatrix}X&0\\0&Y\end{pmatrix}\leq 
\frac{2}{\delta}\begin{pmatrix}I&-I\\-I&I\end{pmatrix}+2\eps\begin{pmatrix}I&0\\0&I\end{pmatrix}+
\mathcal{O}(\kappa),$$ 
for some $\kappa>0$.  In the end we will fix $\sigma$, $\delta$, and
$\eps$ and send $\kappa\ra0$, so we simply ignore the
$\mathcal{O}(\kappa)$-term in the 
following. The first inequality implies 
\begin{align*}
\frac{\sigma m}T\leq
\sup_{i,\alp,\beta}\Big\{&-\tr[\ba(t_0,y_0)Y]+\tr[a(t_0,x_0)X]+\bb(t_0,y_0)p_x-b(t_0,x_0)p_y\\
&+\bc(t_0,y_0)\bu(t_0,y_0)-c(t_0,x_0)u(t_0,x_0)+\bff(t_0,y_0)+f(t_0,x_0)\Big\},
\end{align*}
Note that Lipschitz regularity of the solutions and a standard
argument yields
$$|x_0-y_0|\leq \delta M.$$
So using Ishii's trick on the 2nd order terms
\cite[pp. 33,34]{Is:Unique}, and a few other manipulations, we get
\begin{align*} 
\frac{\sigma m}T\leq
\sup_{i,\alp,\beta}\Big\{&\frac{2}{\delta}|\sigma(t_0,x_0)-\bs(t_0,y_0)|^2 
+2M|b(t_0,x_0)-\bb(t_0,y_0)|\\
&+C\eps(1+|x_0|^2+|y_0|^2)\\
&+M|c(t_0,x_0)-\bc(t_0,y_0)| + |f(t_0,x_0)-\bff(t_0,y_0)|\Big\}.
\end{align*}
Some more work leads to an estimate for $m$ depending on $T$,
$\sigma$, $\delta$, and $\eps$, and using the definition of $m$ and
estimates on $\sup_{i,x,y} \psi_i(0,x,y)$, we obtain a similar upper 
bound for $u-\bu$. We finish the proof of the upper bound on $u-\bu$
by sending $\sigma\ra1$, minimizing this expression w.r.t. $\delta$,
sending $\eps\ra0$, and noting that the result still holds if we
replace $T$ by any $t\in[0,T]$. The lower bound follows in a similar
fashion. 
\end{proof}

\begin{remark}
For more details on such manipulations, we refer to
\cite{JK:ContDep}.
\end{remark}

\end{document}